\documentclass[12pt]{article}

\usepackage[left = 1in, right = 1in, top = 1in, bottom = 1in]{geometry}
\usepackage{amsmath}
\usepackage{amssymb}
\usepackage{latexsym}
\usepackage{bm} 
\usepackage[normalem]{ulem} 
\usepackage{comment}
\usepackage{graphicx}
\DeclareMathOperator*{\argmin}{arg\,min}
\DeclareMathOperator*{\argmax}{arg\,max}
\usepackage{float}
\usepackage{algorithm}
\usepackage{hyperref}
\usepackage{algpseudocode}
\usepackage[dvipsnames]{xcolor}
\usepackage{xspace}

\newcommand{\R}{\mathbb R}
\newcommand{\RV}[1]{{\color{black}#1}}
\newcommand{\eps}{\varepsilon}

\newcommand{\innerprod}[2]{\left\langle #1,#2\right\rangle}
\newcommand{\B}[1]{\bm #1}
\newcommand{\abs}[1]{\left\lvert #1 \right \rvert}

\begin{document}

\begin{center} \large Efficient and Scalable Path-Planning Algorithms for Curvature Constrained Motion in the Hamilton-Jacobi Formulation \end{center}

%% use optional labels to link authors explicitly to addresses:
%% \author[label1,label2]{}
%% \affiliation[label1]{organization={},
%%             addressline={},
%%             city={},
%%             postcode={},
%%             state={},
%%             country={}}
%%
%% \affiliation[label2]{organization={},
%%             addressline={},
%%             city={},
%%             postcode={},
%%             state={},
%%             country={}}

\begin{center} 
Christian Parkinson\footnote{University of Arizona, chparkin@math.arizona.edu (corresponding author)}, Isabelle Boyle\footnote{University of Arizona, ijboyle1@arizona.edu} 
\end{center}
%\author[inst1,inst2]{Author Three}

%\affiliation[UA]{organization={Department of Mathematics, University of Arizona},%Department and Organization
%            addressline={\\617 N. Santa Rita Ave.}, 
%            city={Tucson AZ},
%            postcode={85721}, 
%            \country = {USA}}

\begin{abstract}
%% Text of abstract
We present a partial-differential-equation-based optimal path-planning framework for curvature constrained motion, with application to vehicles in 2- and 3-spatial-dimensions. This formulation relies on optimal control theory, dynamic programming, and Hamilton-Jacobi-Bellman equations. \RV{We develop efficient and scalable algorithms for solutions of high dimensional Hamilton-Jacobi equations which can solve these types of path-planning problems efficiently, even in high dimensions, while maintaining the Hamilton-Jacobi formulation. Because our method is rooted in optimal control theory and has no black box components, it has solid interpretability, and thus averts the tradeoff between interpretability and efficiency for high-dimensional path-planning problems.} We demonstrate our method with several examples. 
\end{abstract}

%%Graphical abstract
%\begin{graphicalabstract}
%\includegraphics{grabs}
%\end{graphicalabstract}

%%Research highlights
%\begin{highlights}
%\item Research highlight 1
%\item Research highlight 2
%\end{highlights}

%\begin{keyword}
%% keywords here, in the form: keyword \sep keyword
%Optimal path-planning \sep Curvature constrained motion \sep Hamilton-Jacobi equations \sep Hopf-Lax formulas \sep Dynamic programming
%% MSC codes here, in the form: \MSC code \sep code
%% or \MSC[2008] code \sep code (2000 is the default)
%\MSC 49L20 \sep 49N90 \sep 93C95 
%\end{keyword}

%% \linenumbers

%% main text
\section{Introduction} \label{Introduction} 

In this manuscript, we develop a Hamilton-Jacobi partial differential equation (PDE) based method for optimal trajectory generation, with special application to so-called Dubins vehicles which exhibit curvature constrained motion. Specifically, we study kinematic models for simple cars, airplanes, and submarines.

Curvature constrained motion was first considered by Dubins \cite{Dubins} who considered a simple vehicle which could only move forward. The model was extended by Reeds and Shepp to a car which could move forward and backward \cite{ReedsShepp}. In both cases, the strategy was to decompose paths into straight line segments and arcs of circles and analyze which combinations could be optimal. Later work in this direction was devoted to adding obstacles \cite{Barra}, and developing algorithms which can produce approximately optimal paths which are robust to perturbation \cite{Agarwal}.

To the authors' knowledge the problem was first analyzed using dynamic programming and PDE by Takei, Tsai and others \cite{TakeiTsai1,TakeiTsai2}. Besides their work, there is a strong precedent in the literature for control theoretic, PDE-based optimal path-planning in a number of applications \cite{Deforest1,OptHumanPathPlanning,OptHumanPathPlanningStochastic,Cars1,CarsMovingObs,Deforest2,Vlad1,Vlad2}. Trajectory generation methods which are rooted in PDE have the advantage that they are easy to implement, entirely interpretable, and can provide theoretical guarantees regarding optimality, robustness, and other concerns. This is to distinguish them from sampling and learning based algorithms (for example \cite{Shukla,Zhang,Wan,Deng,Johnson}) which often sacrifice interpretability for efficiency. The main drawbacks of the PDE-based methods are the lack of efficiency and scalability. Because these methods typically rely on discretizing a domain and approximating a solution to a Hamilton-Jacobi PDE, they can be inefficient even for relatively low-dimensional problems, and intractable for motion planning problems whose state space is more than three dimensions. \RV{While some effort has been made to improve efficiency of grid-based methods through parallelization \cite{Gibou1,Gibou2},} other recent work has been devoted to developing grid-free numerical methods based on Hopf-Lax formulas which can approximate solutions to Hamilton-Jacobi equations efficiently even in high-dimensions \cite{Curse1,Curse2,Lin}. 

\RV{The goal of this manuscript is to develop algorithms for trajectory generation which are (1) interepretable due to the theoretical underpinnings provided by PDEs and optimal control theory, and the absence of any black box components, and (2) efficient enough to be real-time applicable, even in problems with high-dimensional state spaces.} Specifically, we present an optimal path-planning method for simple vehicles which maintains the Hamilton-Jacobi PDE formulation, but does not rely on spatial discretization and is thus efficient and scalable. This requires a reformulation of the standard control theoretic minimal-time path-planning problem analyzed in \cite{TakeiTsai1,TakeiTsai2,Cars1,CarsMovingObs,Deforest2,Vlad1,Vlad2}. \RV{Our method represents a significant step toward fully interpretable PDE-based motion-planning algorithms which are real-time applicable. While the specific focus in this manuscript is Dubins' type vehicles, with small modifications, it is likely that our methods could be adapted to more complex and higher-dimensional models of motion as well.}

The paper is laid out as follows. In section \ref{Modeling}, we present the basic dynamic programming and PDE-based path-planning framework that we use, and introduce models for simple Dubins' type vehicles. In section \ref{Numerics}, we discuss the numerical methods which we use to solve the requisite PDEs and generate optimal trajectories, and their specific application to our problem. In section \ref{results}, we present results of our simulations. We conclude with a brief discussion of our method and potential future research directions in section \ref{Conclusion}. 

\section{Modeling} \label{Modeling}

In this section, we derive the PDE-based path-planning framework which we use. \RV{In particular, because we will solve these PDEs by translating them into optimization problems as described in section \ref{Numerics}, it will be most convenient if we can avoid a formulation which requires boundary conditions, which would translate into difficult constraints in the optimization.} Because of this, we opt for a level-set-type formulation in the vein of \cite{Deforest1,OptHumanPathPlanning,OptHumanPathPlanningStochastic}, as opposed to the control theoretic approach of \cite{TakeiTsai1,TakeiTsai2,Cars1,CarsMovingObs,Deforest2,Vlad1,Vlad2}. These approaches are compatible, but different in philosophy. We compare and contrast them in section \ref{sec:TimeHorizonAndObstacles}.

The level-set method is a general method for modeling contours (or more generally, hypersurfaces of codimension 1) which evolve with prescribed velocity depending on the ambient space and properties inherent to the contour itself \cite{levelSet1,levelSet2}. The basic strategy is to model the contour as the zero level set of an auxiliary function $u: \mathbb R^n \times [0,\infty) \to \mathbb R$ and derive a PDE which $u$ satisfies. As $u$ evolves according to the PDE, the zero level set of $u(\cdot,t)$ evolves, affecting the level set flow. A general, first-order level set equation has the form $$u_t + H(x,\nabla u,t) = 0$$ for some \textit{Hamiltonian} function $H(x,p,t)$ which is homogeneous of degree 1 in the variable $p$.  The level set function $u$ does not need to have any physical meaning, though in many cases (as in ours), level set equations are seen to arise as Hamilton-Jacobi-Bellman equations for feedback control problems where $u$ is a value function. 

We demonstrate this with a brief and formal derivation. Given a time-horizon $T>0$, a starting location $x_0 \in \mathbb R^d$, a desired-ending location $x_f \in \mathbb R^d$, and a function $f: \mathbb R^d\times [0,T] \times \mathbb R^m \to \mathbb R^d$, we consider trajectories $\B x :[0,T] \to \mathbb R^d$ satisfying \begin{equation}  \label{eq:eqOfMotion} \begin{split} &\dot{\bm x} = f(\bm x, t, \bm \alpha), \hspace{1cm} 0 < t \le T, \\ &\bm x(0) = x_0.\end{split}\end{equation} Here $\bm \alpha(\cdot)$ is a control map, taking values in some admissible control set $A \subset \mathbb R^m$. The goal is to choose $\bm \alpha(\cdot)$ so as to steer the trajectory as close as possible to $x_f$ by time $T$. Accordingly, we define the cost functional $$\RV{\mathcal C[\bm x(\cdot)] = \frac 1 2\abs{\bm x(T) - x_f}^2},$$ and we would like to solve the optimization problem $$\inf_{\bm \alpha \in \mathcal A} \mathcal C[\bm x(\cdot)]$$ where $\mathcal A = \{\bm \alpha :[0,T]\to A \, : \, \bm \alpha \text{ measurable}\}$. \RV{Other choices of cost function are possible, but so as to serve as some measure of distance, they should be zero for any path which ends at $x_f$, and increase as $\vert \bm x(T)-x_f\vert$ increases.}

For $x \in \mathbb R^d$ and $t \in [0,T)$, we define the value function \begin{equation}\label{eq:valFuncDef} u(x,t) = \inf_{\bm \alpha \in \mathcal A} \mathcal C_{x,t}[\bm x(\cdot)],\end{equation} where $\mathcal C_{x,t}$ denotes the same functional restricted to trajectories $\bm x(\cdot)$ such that $\bm x(t) = x$. This value function denotes the minimum square distance to the desired endpoint that one can achieve if they are sitting at point $x$ at time $t$. In our case, because there is no running cost along the trajectory, the dynamic programming principle \cite{Bellman} states that the value function is constant along an optimal trajectory. That is, for $\delta > 0$, $$u(\bm x(t),t) = \inf_{\bm \alpha} \{u(\bm x(t+\delta), t+\delta)\},$$ where the infimum is taken with respect to the values of $\bm \alpha(s)$ for $s \in [t,t+\delta)$. If $u$ is smooth, we rearrange, divide by $\delta$ and send $\delta \to 0^+$ to see 
$$\inf_{\alpha \in A} \{u_t + \langle \dot{\bm x},\nabla u\rangle \}= 0.$$  We can use the equation of motion \eqref{eq:eqOfMotion} to replace $\dot{\bm x}$. Further, at time $t = T$, there is no remaining time to travel, so the value is simply the exit cost. Thus we arrive at a terminal-valued Hamilton-Jacobi-Bellman (HJB) equation \begin{equation} \label{eq:HJB}
\begin{split}
&u_t + \inf_{\alpha \in A} \{\langle  f(x,t,\alpha) ,\nabla u\rangle \}= 0,\\
&u(x,T) = \frac 1 2\abs{x-x_f}^2.
\end{split}
\end{equation} The function $H: \R^d \times \R^d \times [0,T] \to \R$ defined by \begin{equation} \label{eq:Ham}
H(x,p,t) = \inf_{\alpha \in A} \{\langle  f(x,t,\alpha) ,p\rangle \}
\end{equation} is the \emph{Hamiltonian}. These computations are formal. There is no reason to believe that $u$ will be smooth, but under very mild conditions on $f$, $u$ will be the Lipschitz continuous viscosity solution of \eqref{eq:HJB} \cite{Viscosity1,Viscosity2,Viscosity3}. Assuming the viscosity solution to \eqref{eq:HJB} is known, one resolves the optimal control map via $$\bm \alpha^*(x,t) = \argmin_{\alpha \in A} \langle f(x,t,a),\nabla u(x,t) \rangle. $$ Again, pending mild assumptions on $f$, this optimization problem has a unique solution whenever $\nabla u(x,t)$ exists, which is true for almost every $(x,t)$ when $u$ is Lipschitz continuous. Finally, one can then resolve the optimal trajectory by integrating \begin{equation}  \label{eq:eqOfMotionOpt} \begin{split} &\dot{\bm x} = f(\bm x, t, \bm \alpha^*(\bm x,t)), \hspace{1cm} 0 < t \le T, \\ &\bm x(0) = x_0,\end{split}\end{equation} or using a semi-Lagrangian method as described in \cite{TakeiTsai2,CarsMovingObs,Vlad2}.

Because we will solve \eqref{eq:HJB} using Hopf-Lax time formulas (and because it is more comfortable for those familiar with PDE), we make the substitution $t \mapsto T-t$, to arrive at an initial value Hamilton-Jacobi-Bellman equation: \begin{equation} \label{eq:HJBIV}
\begin{split}
&u_t - H(x,\nabla u, T-t)= 0,\\
&u(x,0) = \frac 1 2\abs{x-x_f}^2.
\end{split}
\end{equation} For notational convenience, we will still refer to the solution of \eqref{eq:HJBIV} as $u$, though it is a time-inverted version of the solution of \eqref{eq:HJB}. \RV{Also, the Hamiltonians we are interested in below are time-independent (unless moving obstacles are introduced, but these are accounted for separately), so for brevity, we suppress the time-dependence and henceforth write $H = H(x,p)$.}

In the ensuing subsections, we use this basic framework to develop optimal path-planning algorithms for Dubins vehicles in 2- and 3-dimensions. \RV{This 2-dimesional model for kinematic motion goes back to Dubins \cite{Dubins} and Reeds and Shepp \cite{ReedsShepp}. The 3-dimensional versions still serve as simple  models for unmanned airplanes \cite{dubplane1,dubplane2} and seacrafts \cite{dubsub2} today. We note that while our inspiration is  self-driving vehicles, curvature constrained motion has other applications such as maneuvering of bevel tipped needles through biological tissue \cite{needles}.}

\subsection{Dubins Car}

 We consider a simple rectangular car as pictured in figure \ref{fig:car}. We let $(x,y)\in\mathbb R^2$ denote the center of mass of the car and $\theta \in [0, 2\pi]$ denote the orientation, measured counterclockwise from the horizontal. We refer to $(x,y,\theta)$ as the \textit{configuration} of the car. The car has a rear axle of length $2R$ and distance $d$ from the center of mass to the center of the rear axle. The motion of the car is subject to the nonholonomic constraint
\begin{equation}
\dot{\bm y}\cos{\bm \theta} - \dot{\bm x}\sin{\bm \theta} = d\dot{\bm \theta} \label{eq:HolonomicConstraint}
\end{equation} which specifies that movement occurs tangential to the rear wheels. Motion is also constrained by a maximum angular velocity  $|\dot{\theta}| \leq W$ (or equivalently a minimum turning radius $\rho = 1/W$). The kinematics for the car are given by:

\begin{figure}[t!]
\centering
\includegraphics[width = 0.5\textwidth]{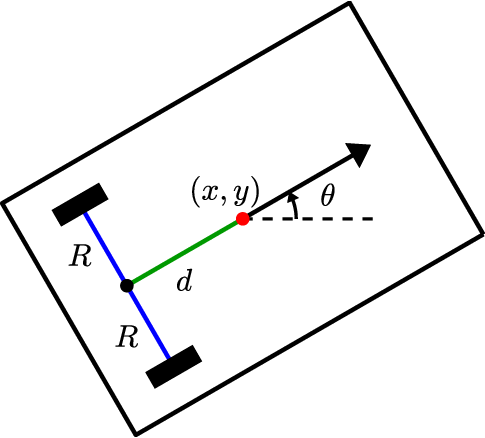}
\caption{A simple self-driving car}
\label{fig:car}
\end{figure}

\begin{align}
\begin{split}
&\dot{\bm{x}} = \bm{v}\cos{\B \theta} - \bm{\omega}Wd\sin{\bm{\theta}}, \\
&\dot{\bm{y}} = \bm{v}\sin{\B \theta} + \bm{\omega}Wd\cos{\bm{\theta}}, \\
&\dot{\bm{\theta}} = \bm{\omega}W.
\end{split}
\label{eq:kinematicEquationsCar}
\end{align}
Here, $\bm{v}(\cdot), \bm{\omega}(\cdot) \in [-1,1]$ are  normalized control variables representing tangential and angular velocity, respectively. 

Given a desired final configuration $(x_f,y_f,\theta_f)$ for the car, we can insert these dynamics into the above derivation so that the (time-reversed) optimal remaining distance function $u(x,y,\theta,t)$ for this problem satisfies the Hamilton-Jacobi-Bellman equation 
$$u_t - \inf_{v,\omega} \Big\{u_x(v\cos\theta - \omega Wd \sin \theta) + u_y(v\sin \theta + \omega W d \cos \theta) + \omega W u_\theta \Big\} = 0.$$ Rearranging to isolate $v,\omega$, we have $$u_t - \inf_{v,\omega} \Big\{v(u_x \cos \theta + u_y \sin \theta) + \omega Wd(-u_x \sin \theta + u_y \cos \theta + \tfrac 1d u_\theta) \Big\} = 0.$$ We recall that the control variables $v,\omega$ are normalized: $v,\omega \in [-1,1]$. Thus since the above infimum is linear in $v$ and $\omega$ and their dependence is decoupled, the infimum for each control variable must occur at one of the endpoints. Assuming $u(x,y,\theta,t)$ is known, the optimal controls are given by\begin{equation}\label{eq:carControls}
\begin{split}
v^*(x,y,\theta,t) &= -\text{sign}(u_x(x,y,\theta,t)\cos\theta + u_y(x,y,\theta,t)\sin\theta),\\
\omega^*(x,y,\theta,t) &= -\text{sign}(-u_x(x,y,\theta,t)\sin\theta + u_y(x,y,\theta,t)\cos\theta + \tfrac 1 d u_\theta).
\end{split}
\end{equation}
and 
\begin{equation} \label{eq:HJBCar1}
\begin{split}
&u_t + |u_x\cos{\theta} + u_y\sin{\theta}| + Wd|-u_x\sin{\theta} + u_y\cos{\theta} + \tfrac 1 d u_{\theta}| = 0, \\
&u(x,y,\theta,0) = \frac 1 2|(x,y,\theta) - (x_f,y_f,\theta_f)|^2.
\end{split}
\end{equation}
\RV{As written, one can specify a final location \emph{and} orientation that the car must achieve (hence the inclusion of $\theta_f$ in the initial condition). If one only cares to specify a final location, this can be omitted, and the initial condition can read $u(x,y,\theta,0) = \frac 1 2 \abs{(x,y)-(x_f,y_f)}^2$. A similar observation applies for the 3-dimensional models below.} When modeling the car as a point mass (i.e. setting $d=0$), \eqref{eq:HJBCar1} simplifies to
\begin{equation}\label{eq:HJBCar}
u_t + |u_x\cos{\theta} + u_y\sin{\theta}| + W|u_{\theta}| = 0.
\end{equation} 

We will deal with this latter formulation (i.e., $d= 0$) so that we are neglecting the actual shape of the vehicle as in \cite{TakeiTsai1,TakeiTsai2}. Note that, as presented, we have not yet accounted for obstacles. It is argued in \cite{Cars1,CarsMovingObs} that it is easier to account for obstacles when one does not simplify the car to a point mass, but the formulation is slightly different in those papers, and we will need to make special considerations for obstacles regardless. We discuss this further in section \ref{sec:TimeHorizonAndObstacles}. In the next two subsections, we present two higher dimensional generalizations of this model. 

\subsection{Dubins Airplane}

We generalize to the case of a simple airplane by introducing a third spatial variable to the dynamics. As a result the autonomous vehicle operates in the coordinates $(x,y,z,\theta)$, rather than just $(x,y,\theta)$. This is achieved by adding in a fourth kinematic equation, which applies a restriction to the motion in the $z$ direction akin to that imposed on $\theta$. For the sake of this model, the constraint imposed on $z$ is decoupled from the constraint on the planar motion, in a manner similar to that of an airplane so we call this the \emph{Dubins Airplane} model, though, as in the case of the car, we are simplifying the dynamics compared to that of a real airplane. We once again consider the vehicle to be a point mass. This leads to the kinematic equations 
\begin{align}
\begin{split}
&\dot{\bm{x}} = \bm{v}\cos{\B \theta}, \\
&\dot{\bm{y}} = \bm{v}\sin{\B \theta}, \\
&\dot{\bm{z}} = \bm{\omega_z}W_z, \\
&\dot{\bm{\theta}} = \bm{\omega_{xy}}W_{xy}, 
\end{split}
\label{eq:kinematicEquationsPlane}
\end{align}
where, once again, we enforce  $\bm{v}(\cdot), \bm{\omega}_{xy}(\cdot), \bm{\omega}_{z}(\cdot) \in [-1,1]$. In \eqref{eq:kinematicEquationsPlane}, $W_{xy} > 0$ is a constraint on the angular velocity in the $xy$-plane, and $W_z > 0$ is a constraint on the angular velocity in the vertical direction. 

Similar reasoning will lead us to the Hamilton-Jacobi-Bellman equation
\begin{equation}
\begin{split}
&u_t + |u_x\cos{\theta} + u_y\sin{\theta}| + W_z|u_{z}|+ W_{xy}|u_{\theta}|  = 0, \\
&u(x,y,z,\theta,0) = \frac 1 2|(x,y,z,\theta) - (x_f,y_f,z_f,\theta_f)|^2.
\label{eq:HJBPlane1}
\end{split}
\end{equation} To truly model an airplane in flight, it makes sense to restrict to unidirectional velocity: $\B v(\cdot) =1$. In this case, the equation simplifies to \begin{equation} \label{eq:HJBPlane}
u_t - u_x\cos{\theta} - u_y\sin{\theta} + W_z|u_{z}|+ W_{xy}|u_{\theta}| = 0.
\end{equation}

As mentioned above, this formulation does not yet account for the presence of obstacles, which will be discussed in section \ref{sec:TimeHorizonAndObstacles}. One could also account for finer-scale modeling concerns such as minimum cruising velocity, but for our purposes, the formulation is as presented.

\subsection{Dubins Submarine}

As an alternative to Dubins Airplane, we can model the problem in three dimensions by enforcing a total curvature constraint, rather than decoupling the constraints on the planar and vertical motion. We dub this model the \emph{Dubins Submarine}, though again, we are neglecting many of the dynamics of a real submarine. In this case, we let $(x,y,z)$ represent the center of mass of the vehicle. The orientation of the vehicle is now represented by a pair of angles: $\theta \in [0,2\pi]$ represents the angular orientation on the $xy$-plane (the azimuthal angle), and $\phi \in [0,\pi]$ is the angle of inclination, with $\phi = 0$ pointing straight up the $z$-axis and $\phi = \pi$ pointing straight down the $z$-axis. In this case, the equations of motion are \begin{align}
\begin{split}
&\dot{\bm{x}} = \bm{v}\cos{\theta}\sin\varphi, \\
&\dot{\bm{y}} = \bm{v}\sin{\theta}\sin\varphi, \\
&\dot{\bm{z}} = \bm{v}\cos\varphi,\\
&\dot{\bm{\theta}} = \bm{\omega}_1W, \\
&\dot{\bm{\varphi}} = \bm{\omega}_2W,
\end{split}
\label{eq:kinematicEquationsSub}
\end{align}
where $\B v(\cdot)$ is the normalized control variable representing tangential velocity, $\B \omega_1(\cdot),\B \omega_2(\cdot)$ are normalized control variables for angular velocity, and $W >0$ is a constraint on the curvature of the path. The magnitude of the curvature of a path obeying \eqref{eq:kinematicEquationsSub} is given by $$\kappa(t) := \left\vert \frac{d}{dt}\frac{(\dot{\bm{x}},\dot{\bm{y}},\dot{\bm{z}})}{\vert(\dot{\bm{x}},\dot{\bm{y}},\dot{\bm{z}})\vert} \right\vert= W\sqrt{\bm{\omega}_1^2 \sin^2(\B \varphi) + \B \omega_2^2},$$ which leads to a contraint on $\B \omega_1, \B \omega_2$ of the form \begin{equation} \label{eq:curvConstraint}\sqrt{\bm{\omega}_1^2 \sin^2(\B \varphi) + \B \omega_2^2}\le 1.\end{equation}

Going through the same derivation for the Hamilton-Jacobi-Bellman equation (and having made the substitution $t\mapsto T-t$ to reverse time), we arrive at $$u_t - \inf_{v,\omega_1,\omega_2} \left\{ v(u_x \cos\theta \sin\varphi + u_y \sin \theta \sin \varphi + u_z \cos \varphi)+ \omega_1 u_\theta + \omega_2 u_\varphi \right\} = 0.$$ Again, the minimization in $v$ is decoupled from that of $(\omega_1,\omega_2)$, so it is resolved exactly as before. For the minimization in $(\omega_1,\omega_2)$, assuming that $\sin\varphi \neq 0$, we write $$\inf_{\omega_1,\omega_2}\{\omega_1 u_\theta + \omega_2 u_\varphi\} = \inf_{\omega_1,\omega_2} \left\{ (\omega_1\sin\varphi)\left(\frac{u_\theta}{\sin \varphi}\right) + \omega_2 u_\varphi\right\}.$$ Recall, $(\omega_1,\omega_2)$ are constrained by \eqref{eq:curvConstraint} so that $(\omega_1\sin \varphi,\omega_2)$ is in the unit circle. Thus, $$\inf_{\omega_1,\omega_2} \left\{ (\omega_1\sin\varphi)\left(\frac{u_\theta}{\sin \varphi}\right) + \omega_2 u_\varphi\right\} = -\sqrt{\frac{u_\theta^2}{\sin^2 \varphi} + u_\varphi^2}.$$ Thus the HJB equation for this model is
\begin{equation}
\begin{split}
&u_t + |u_x\cos{\theta}\sin{\varphi} + u_y\sin{\theta}\sin{\varphi} + u_z\cos{\varphi}| + W\sqrt{\frac{u_\theta^2}{\sin^2{\varphi}} + u_\varphi^2 } = 0, \\
&u(x,y,z,\theta,\varphi,0) = \frac 1 2|(x,y,z,\theta,\varphi) - (x_f,y_f,z_f,\theta_f,\varphi_f)|^2,
\label{eq:HJBSub}
\end{split}
\end{equation} for $(x,y,z,\theta,\varphi,t) \in \mathbb R^3 \times [0,2\pi]\times (0,\pi) \times [0,T].$ When $\varphi = 0$ or $\varphi = \pi$, the submarine is oriented directly upward or downward respectively. In this case, to change the orientation, one must modify $\phi$ (i.e., $\omega_1$ can take any value, but it will not affect the orientation; only $\omega_2$ affects the orientation). One could specifically account for this in the above derivation if desired. For computational purposes, it suffices to replace the $\sin^2 \varphi$ in the denominator in \eqref{eq:HJBSub} with $(\sin^2 \varphi + \varepsilon)$ where $\varepsilon$ is only a few orders of magnitude larger than machine precision. This will circuit any division by zero without materially affecting results. For our purposes, we use $\varepsilon = 10^{-10}.$

In the ensuing sections, we discuss the formulation of these models in the presences of impassable obstacles, and then move on to develop numerical methods for approximating these equations and generating optimal trajectories. 

\subsection{Level Set vs. Optimal Control Formulation: the Time Horizon and Obstacles} \label{sec:TimeHorizonAndObstacles}

In all the preceding derivations, we implicitly assume that the vehicles are moving in free space. Here we discuss how one may account for obstacles which impede the vehicles. We also discuss the role of the time horizon $T$, because both the manner in which we include obstacles and the time horizon arise as a consequence of modeling the problem using level set equations, as opposed to what is perhaps a more natural control-theoretic model. \RV{To reiterate, we use level set equations because they are amenable to the available numerical methods. Specifically, using the level set formulaton allows us to derive an HJB equation which can resolve time-optimal paths \emph{without} using any boundary conditions. This is desirable because in section \ref{Numerics}, we describe how one can solve the HJB equation by reframing it as an optimization problem for which it difficult to incorporate boundary conditions. However, this decision has ramifications that deserve some discussion.}

Reverting back to the general framework from the beginning of section \ref{Modeling}, we recall that for our models, the value function $u(x,t)$ encodes the minimum achievable distance to the desired endpoint given that the vehicle is in position $x$ at time $t$. An alternative approach to path-planning, like that used in \cite{TakeiTsai1,TakeiTsai2,Cars1,CarsMovingObs,Deforest2,Vlad1,Vlad2}, is to define the value function $\tau(x,t)$ to be the optimal remaining travel time to the desired endpoint given that the vehicle is in position $x$ at time $t$. In the latter case, by definition, $\tau(x,t) = +\infty$ whenever there is no admissible path for a car which is at $x$ at time $t$ which can steer the car to the desired ending point before hitting the time horizon $T$, and thus $\tau(x,t)$ is finite if and only if there is an admissible path to the ending point in the allotted time. When there are no obstacles (or stationary obstacles), once $\tau(x,t)$ becomes finite, it will remain constant. For example, if the optimal travel time from $x$ to $x_f$ is $3$ units of time, then $\tau(x,1) = +\infty$, $\tau(x,3) = 3$ and $\tau(x,10) = 3$. In this case, one can essentially eliminate the time horizon, by simply taking $T$ large enough that $\tau(x,T) < +\infty$ for all points $x$ that one cares about. Specifically, with mild assumptions on the dynamics (so that admissible paths are possible from any given point), for any compact spatial domain $\Omega$, there is a time horizon $T^* = T^*(\Omega) > 0$, such that $\tau(x,T^*) < +\infty$ for all $x \in \Omega,$ and choosing any $T \ge T^*$ as the time horizon for the problem will yield the same results. In this way, it does not matter that one optimally chooses the time horizon: as long as it is large enough, one can resolve the time-optimal path from $x$ to $x_f$ and this path will be independent of the time horizon $T$. Moreover, with this modeling choice, whenever $\tau(x,T)$ is differentiable, one can uniquely determine the optimal control values and thus there is a unique time optimal path from $x$ to $x_f$. 

The interpretation of our model is slightly different since $u(x,t)$ denotes some measure of optimal distance to the endpoint. Like $\tau(x,t)$, given stationary obstacles and mild conditions on the dynamics, this function will become constant in finite time: for any compact spatial domain $\Omega,$ there is a time horizon $T^* = T^*(\Omega) >0$ such that $u(x,T^*) = 0,$ and then $u(x,t) = 0$ for all $t > T^*.$ In essence, this says that given enough time, there are paths which can reach the desired final point. However, in this case, optimality of a trajectory is judged solely in view of ``distance to the desired endpoint" so any path that reaches the final endpoint is equally optimal. Because of this, choosing the time horizon optimally becomes important. Specifically, given a point $x$, define $t_x = \inf\{t > 0 \, : \, u(x,t) = 0\}.$ If $T > t_x$ (and if the vehicle satisfies a small-time local controllability condition \cite{STLC}), then there will be infinitely many ``optimal" paths from $x$ to $x_f$ since there is more time that needed in order to reach the final point. In this case, we would like the optimal path which requires the minimal time to traverse. In order to resolve the minimal time path beginning from a point $x$, we need to actually use $t_x$ as the time horizon. 

This causes some difficulty because resolving $t_x$ for a given point $x$ requires reformulating the problem in terms of $\tau(x,t)$ as discussed above. Empirically, this affects the numerical methods as well, as we discuss in section \ref{results}. Before this, we address one further modeling concern.

One last modeling concern is the inclusion of obstacles. Intuitively, any path the intersects with an obstacle at any time should be considered illegal and assigned infinite cost, so that it is never the optimal path. To model this, assume that a vehicle is navigating a domain $\Omega \subset \R^n$ which, at any time $t >0$ is disjointly segmented into free space and obstacles, $\Omega = \Omega_{\text{free}}(t) \sqcup \Omega_{\text{obs}}(t).$  

Using the ``travel-time" formulation of \cite{TakeiTsai1,TakeiTsai2,Cars1,CarsMovingObs} as discussed above, one assigns infinite cost to paths that intersect with obstacles by setting $\tau(x,t) = +\infty$ for any points $(x,t) \in \Omega \times (0,T]$  such that $x \in \Omega_{\text{obs}}(t)$. This becomes a crucial boundary condition in the HJB equation for the travel time function $\tau.$ However, as described in section \ref{Numerics}, we hope to resolve the solution to our HJB equations using optimization routines which are much easier to implement when the HJB equations are free of boundary conditions. This is one of the primary motivations for using the level set formulation. In the level set formulation, one incorporates obstacles not by enforcing a boundary condition, but by setting velocity to zero in the obstacles. That is, define $O(x,t)$ to be the indicator function of the free space at time $t$: \begin{equation} \label{eq:OBS}
O(x,t) = \bm{1}_{\Omega_{\text{free}}(t)}(x) = \begin{cases}0, & x \in \Omega_{\text{obs}}(t), \\ 1, & x \in \Omega_{\text{free}}(t),  \end{cases} \,\,\,\,\,\,\, (x,t) \in \R^n \times [0,T].
\end{equation} To set velocity to zero in the obstacles, we then simply multiply the Hamiltonians in \eqref{eq:HJBCar},\eqref{eq:HJBPlane}, \eqref{eq:HJBSub} (or generally in \eqref{eq:HJBIV}) by $O(x,t)$. For example, the version of \eqref{eq:HJBCar} that we actually solve is \begin{equation} \label{eq:HJBCarObs}
u_t + O(x,y,t)\big[|u_x\cos{\theta} + u_y\sin{\theta}| + W|u_{\theta}|\big] = 0.
\end{equation} \RV{Having done this, any path which enters an obstacle will be forced to stop, get stuck, and fail to reach the end point, which will never be optimal.} In this way, we have incorporated obstacles without using any boundary conditions. This raises other numerical issues (for example, how to efficiently determine whether a given $x$ lies in a obstacle at time $t$) which we discuss further in the ensuing section. \RV{This manner of incorporating obstacles is akin to that of \cite{TakeiTsai2, Deforest1, OptHumanPathPlanning}.}

\section{Numerical Methods} \label{Numerics}

To this point, the standard approach to approximating solutions to HJB equations in applications like this has been to use finite difference schemes such as fast-sweeping schemes \cite{fastSweep1,fastSweep2,fastSweep3}, fast-marching schemes \cite{FastMarch1,FastMarch2}, and their generalizations \cite{FastMarch3,OUM}. These are used, for example, in \cite{TakeiTsai1, TakeiTsai2, OptHumanPathPlanning, Cars1, CarsMovingObs}. These methods are easy to implement and can be adapted for high-order accuracy, but because they are grid-based, their computational complexity scales exponentially with the dimension of the domain, meaning that they are only feasible in low dimension. \RV{The authors of \cite{Gibou1,Gibou2} present parallelizable fast-sweeping schemes which can approximate the solutions of 3- and 4-dimensional HJB-type equations in a matter of minutes (or in the case of the Eikonal equation, a matter of seconds).    However, because the model for the Dubins' submarine is 5-dimensional, and we would like to develop methods which scale to even higher dimension, we opt for a grid-free scheme.}

Recent numerical methods attempt to break the \emph{curse of dimensionality} \cite{Bellman} by resolving the solution to certain Hamilton-Jacobi equations at individual points using variational Hopf-Lax type formulas \cite{Curse1,Curse2}. \RV{In particular, given a Hamiltonian $H: \mathbb R^d \times \mathbb R^d \to \R$ which is convex in its second argument and $C^2$-smooth, and an initial data function $g:\R^d \to \R$ which is convex and coercive,} the authors of \cite{Curse2} conjecture that the solution to the state-dependent Hamilton-Jacobi equation \begin{equation} \label{eq:HJgen} u_t + H(x,\nabla u) = 0, \,\,\,\,\,\,\,\,\, u(x,0) = g(x),\end{equation} is given by the generalized Hopf-Lax formula, 
\begin{equation} \label{eq:minProb}
\begin{split}
u(x,t) &= \min_{v\in\mathbb{R}^d}\left\{g(\bm{x}(0)) + \int_0^t \langle \bm{p}(s), \nabla_p H(\bm{x}(s),\bm{p}(s)) -  H(\bm{x}(s),\bm{p}(s))  \rangle ds\right\}\\
&\hspace{-5mm}\text{subject to the Hamiltonian dynamics,}\\
&\dot{\bm{x}}(s) = \hphantom{-}\nabla_p H(\bm{x}(s),\bm{p}(s)), \,\,\,\,\,\,\,\,\,\,\, \bm{x}(t) = x, \\
&\dot{\bm{p}}(s) = -\nabla_x H(\bm{x}(s),\bm{p}(s)), \,\,\,\,\,\,\,\,\,\,\, \bm{p}(t) = v.
\end{split}
\end{equation} \RV{While the conjecture requires the Hamiltonian to be $C^2$, the authors of \cite{Curse2} give specific examples which violate this assumption for which the conjecture still holds, so it may hold even more broadly.}

\RV{Similarly, beginning from \eqref{eq:valFuncDef}, we can derive a conjectured saddle-point formula to approximate the value function. To do so, we discretize the time interval $t=t_0<t_1<\cdots < t_N = T$ and let $x_j$ be discrete approximations to $\bm x(t_j)$ for $j=0,1,\ldots, N$. For simplicity, we will always assume a uniform discretization $t_j = t + j\delta$, $j=0,1,\ldots,N$ where $\delta = \frac{T-t}{N}$, though this is not strictly necessary. Then we see \begin{equation} \label{eq:saddlePoint1} 
\begin{split} 
u(x,t) &= \inf_{\bm x(\cdot),\bm \alpha(\cdot)} \left\{g(x(T)) \, : \, \bm x(t) = x, \dot{\bm x}=f(\bm x,s,\bm\alpha), s \in (t,T]  \right\}, \\ 
\end{split}
\end{equation} where $g(x) = \frac 12\abs{x-x_f}^2$. Formally, discretizing the constraint using a backward Euler scheme, letting $\alpha_j = \bm\alpha(t_j)$ and $f_j = f(x_j,t_j,\alpha_j)$, and introducing Lagrange multipliers $p_j$, we have \begin{equation}\label{eq:saddlePoint2} \begin{split} u(x,t) &\approx \inf_{x_j,\alpha_j} \left\{g(x_N) \, : \, x_0 = x, x_{j} = x_{j-1} + \delta f_j, \,\,\, 1\le j \le N\right\} \\
&=\inf_{x_j,\alpha_j}\sup_{p_j} \left\{g(x_N) + \sum^{N}_{j=1} \innerprod{p_j}{\delta f_j + x_{j-1}- x_{j}}\right\} \\ 
&=\inf_{x_j,\alpha_j}\sup_{p_j} \left\{g(x_N) +  \sum^{N}_{j=1} \innerprod{p_j}{x_{j-1}-x_{j}} + \delta \sum^{N-1}_{j=0} \innerprod{p_j}{f_j}\right\}.
\end{split} \end{equation} From here, we interchange the infimum over $\alpha_j$ with the supremum over $p_j$ and write \begin{equation}
\label{eq:saddlePoint2}
\begin{split}
u(x,t) &\approx \inf_{x_j}\sup_{p_j} \left\{g(x_N) +  \sum^{N}_{j=1} \innerprod{p_j}{x_{j-1}-x_{j}} + \delta \sum^{N-1}_{j=0} \inf_{\alpha_j} \{\innerprod{p_j}{f_j}\}\right\} \\
&= \inf_{x_j}\sup_{p_j}\left\{g(x_N) +  \sum^{N}_{j=1} \innerprod{p_j}{x_{j-1}-x_{j}} + \delta \sum^{N-1}_{j=0} H(x_j,p_j)\}\right\} 
\end{split}
\end{equation} Finally, as in \eqref{eq:HJBIV}, making the time-reversing substitution $t\mapsto T-t$ yields  \begin{equation}
\label{eq:saddlePoint3}
u(x,t) \approx\inf_{x_j}\sup_{p_j}\left\{g(x_0) +  \sum^{N}_{j=1} \innerprod{p_j}{x_{j}-x_{j-1}} - \delta \sum^{N-1}_{j=0} H(x_j,p_j)\}\right\}.
\end{equation}
In general, swapping the order of optimization does not yield the same value. However, performing this step provides a connection between this optimization problem and the Hamilton-Jacobi-Bellman equation via the Hamiltonian. In the spirit of the conjectured Hopf-Lax formula \eqref{eq:minProb}, the authors of \cite{Lin} conjecture (and provide solid empirical evidence) that having swapped the order of the optimization, the resulting saddle-point problem still provides an approximation of the soution of the HJB equation. Our results corroborate this. In fact, \cite{Lin} specifically considers Eikonal-type equations whose Hamiltonians $H(x,p) = V(x)\abs{p}$ share key properties with ours, such as convexity and 1-homogeneity in the second variable and smoothness almost everywhere, so there is precedent to expect this saddle point problem can indeed represent our solutions.}

An alternating minimization technique in the spirit of the primal-dual hybrid gradient method \cite{PDHG1,PDHG2,PDHG3} is proposed by \cite{Lin} to solve this saddle-point problem. For completeness of our exposition, we reprint this in Algorithm \ref{alg:1} (note: this is a reprint of Algorithm 5 in \cite{Lin}, modified slightly to fit our equations). As this algorithm resolves the values of $u$ (the solution of \eqref{eq:HJgen}), it also resolves discrete approximations of the optimal trajectory $\bm x(s)$, and the optimal costate trajectory $\bm p(s)$, which can be seen as a proxy for $\nabla u$ along the optimal trajectory. 

\begin{algorithm}[t!]
\caption{Splitting Method for Solving \eqref{eq:minProb}}
\hspace*{\algorithmicindent} Given a point $(x,t) \in \mathbb R^d \times (0,T]$, a Hamiltonian $H$, an initial data function $g$, a time-discretization count $N$, a max iteration count $K$, an error tolerance TOL, and relaxation parameters $\sigma,\tau,\kappa > 0$, we resolve the minimization problem \eqref{eq:minProb} as follows. \\

Set $x^1_N = x, p^1_0 = 0$ and $\delta = t/N$. Initialize $x^1_0,x^1_1,\ldots,x^1_{N-1},p^1_1,\ldots, p^1_N$ randomly, and set $z^1_j = x^1_j$ for all $j = 0,1,\ldots, N$. 
\begin{algorithmic}
\For {$k = 1$ to $K$}
    \State Set $p^{k+1}_0 = 0$
    \For {$j = 1$ to $N$}
    \State $p_j^{k+1} = \argmin_{\tilde p} \{ \delta \sigma H(x_j^k,\tilde p) + \frac{1}{2} \left\vert \tilde p-(p_j^k + \sigma(z_j^k - z_{j-1}^k)) \right\vert_2^2 \}$
    \EndFor
    \State $x_0^{k+1} = \argmin_{\tilde x} \{ \tau g(\tilde x) + \frac{1}{2} \left\vert \tilde x - (x_0^k + \tau p_1^{k+1}) \right\vert_2^2 \}$ \,\,\,\,\,\,\,\,\, (note: $p^{k+1}_0 = 0$) 
    \For{$j = 1$ to $N - 1$}
    \State $x_j^{k+1} = \argmin_{\tilde x}\{ - \delta\tau H(\tilde x,p_j^{k+1}) + \frac{1}{2} \left\vert \tilde x - (x_j^k - \tau(p_j^{k+1} - p_{j+1}^{k+1})) \right\vert_2^2 \}$
    \EndFor
    \State Set $x^{k+1}_N = x$ 
    \For{$j = 0$ to $N$}
    \State $z_j^{k+1} = x_j^{k+1} + \kappa(x_j^{k+1} - x_j^k)$
    \EndFor
    \State change $= \max \{ \lVert x^{k+1} - x^k \rVert , \lVert p^{k+1} - p^k \rVert \} $
    \If {change $<$ TOL}
    \State break
    \EndIf
\EndFor
\State $u = g(x_0) + \sum_{j=1}^N \langle p_j, x_j - x_{j-1} \rangle - \delta H(x_j,p_j)$
\State \textbf{return } $u$; the value of the solution of \eqref{eq:HJgen} at the point $(x,t)$
\end{algorithmic}
\label{alg:1}
\end{algorithm}

\RV{Assuming a saddle point exists and that the minima are resolved exactly, convergence of the algorithm to a saddle point is guaranteed when $\sigma\tau \le 0.25$ and $\kappa \in [0,1]$ \cite{Lin,PDHG1}. In our particular implementation, we take $\sigma = 0.5$, $\tau = 0.5$, and $\kappa = 1$. } The norm used to determine convergence is not terribly important. If the state space is $d$-dimensional, so that at each time-step $j = 0,1\ldots,N$ and iteration $k$, we have $x^k_j = ((x_1)_j^k, (x_2)_j^k,\ldots,(x_{d})_j^k)$, we use the norm $$\|x^{k+1}- x^k\| = \sup_{1\le i \le d, 0 \le j \le N} \abs{(x_i)^{k+1}_j - (x_i)^k_j}$$ and similarly for $p$. Thus we halt the iteration when no coordinate of either $x$ or $p$ has changed by more than the prescribed $\text{TOL}$ value. 

Note that at each iteration in algorithm \ref{alg:1}, there are approximately $2N$ optimization problems \begin{align}
    p_{j}^{k+1}&=\argmin_{\tilde p} \{ \delta \sigma H(x_j^k,\tilde p) + \frac{1}{2} \left\vert \tilde p-(p_j^k + \sigma(z_j^k - z_{j-1}^k)) \right\vert_2^2 \}, \label{eq:pmincar}\\
    x_0^{k+1}&=\argmin_{\tilde x} \{ \tau g(\tilde x) + \frac{1}{2} \left\vert \tilde x - (x_0^k + \tau p_1^k)) \right\vert_2^2 \}, \label{eq:x0mincar}\\
    x_j^{k+1}&=\argmin_{\tilde x}\{ -\delta \tau H(\tilde x,p_j^{k+1}) + \frac{1}{2} \left\vert \tilde x - (x_j^k - \tau(p_j^{k+1} - p_{j+1}^{k+1})) \right\vert_2^2 \},\label{eq:xmincar}
\end{align} where $N$ is the number of discrete time steps along the path. These are vector valued quantities, with the subscript denoting the time step along the path and the superscript denoting the iteration number. We allow a maximum of $K =100000$ in our simulations, though with baseline parameter values, the routine usually converged to within a tolerance of $\text{TOL} = 10^{-3}$ within 10000 iterations and never reached the maximal iteration count (we discuss this more specifically in section \ref{results}). Even so, this is a fairly large computational burden, so when possible, it behooves one to resolve the optimization problems in algorithm \ref{alg:1} analytically. When this is not possible, one can use gradient descent or some other comparably simple optimization method. In our applications, the minimization problems in the costate variables $p$ can be resolved exactly, while the minimization problems in the state variables $x$ will need to be approximated in some cases. Empirically, it was observed in \cite{Lin} (and corroborated in our simulations) that one only needs to very crudely approximate these minimizers; for example, using only a few steps of gradient descent at each iteration. In this manner, the approximation may be very poor at early iterations, but becomes better as the iteration count increases.  We describe the specifics of how we resolve the minimization problems from algorithm \ref{alg:1} in the next subsections.

\subsection{Implementation of algorithm \ref{alg:1} for the Dubins car}\label{car_updates}

Here we describe the implementation of algorithm \ref{alg:1} for the  Dubins car. In \eqref{eq:HJBCar}, we see that the Hamiltonian is
\begin{equation} \label{eq:carHam} H(x,p) = |p_1 \cos(x_3) + p_2 \sin(x_3)| + W|p_3|,\end{equation} where $ (x_1,x_2,x_3)$ represent $(x,y,\theta)$ respectively and $(p_1,p_2,p_3)$ are proxies for $(u_x,u_y,u_\theta)$ respectively. 

An iteration in algorithm \ref{alg:1} begins by resolving the minimization problem \eqref{eq:pmincar}, which can be done analytically. For completeness, we include the formal derivation of the minimizer here. For this Hamiltonian, the dependence of $H$ on $(p_1,p_2)$ is decoupled from the dependence on $p_3$, so we can treat these as two separate problems. In what follows, all indexing notation is adapted to MATLAB indexing conventions so that $\tilde p_{1:2}$ indicates the first two components of the vector $\tilde p$, for example. 

To simplify notation, we set $$\gamma = (\cos((x_3)_{j}^k), \sin((x_3)^k_{j})), \,\,\,\,\,\, \beta = p_j^{k} + \sigma \left(z_j^{k} - z_{j-1}^{k} \right),$$ so that \eqref{eq:pmincar} is can be written $$ p^{k+1}_j = \argmin_{\tilde p}\left\{\delta \sigma  \abs{\gamma^T\tilde p_{1:2}} + W\abs{\tilde p_3} + \frac{1}{2} \| \tilde p - \beta\|^2_2\right\}.$$ Note that, when resolving $p^{k+1}_{j}$, $\gamma$ and $\beta$ are known quantities which depend on $j$ and $k$, but we suppress this dependence for ease of notation. Taking the gradient of the function being minimized and setting to zero, we see that the minimizer $\tilde p$ satisfies
\begin{equation}\label{eq:pDerivation}
\delta \sigma \frac{\gamma^T \tilde p_{1:2} }{\abs{\gamma^T\tilde p_{1:2}}} \gamma + \tilde p_{1:2} - \beta_{1:2} = 0.
\end{equation}
so that 
\begin{equation*}
\left( \frac{\delta \sigma}{\abs{\gamma^T\tilde p_{1:2}}}\gamma\gamma^T + I \right)\tilde p_{1:2} = \beta_{1:2}.
\end{equation*}
Now in order to solve for $\tilde p_{1:2}$, we project all vectors along $\gamma$ and the orthogonal direction to $\gamma$ by writing,
\begin{equation}\label{eq:pDerivative2}
\tilde p_{1:2} = r\gamma + s(\beta_{1:2} - (\gamma^T \beta_{1:2})\gamma),
\end{equation} for some $r,s \in \R$, and  
\begin{equation*}
\beta_{1:2} = (\gamma^T\beta_{1:2}) \gamma + (\beta_{1:2} - (\gamma^T \beta_{1:2} )\gamma).
\end{equation*}
Inserting these in \eqref{eq:pDerivation} and using $\gamma^T \gamma = 1$ gives
\begin{equation*}
\frac{\delta \sigma r}{\abs{r}}\gamma + r\gamma + s(\beta_{1:2} - (\gamma^T \beta_{1:2}) \gamma) = (\gamma^T\beta_{1:2}) \gamma + (\beta_{1:2} - (\gamma^T\beta_{1:2} )\gamma).
\end{equation*}
whereupon we immediately have $s = 1$. Then
\begin{equation*}
\left( \frac{\delta \sigma}{\abs{r}} + 1 \right)r = \gamma^T\beta_{1:2}
\end{equation*} shows that $r$ shares a sign with $\gamma^T\beta_{1:2}$, so we can write $r = a (\gamma^T\beta_{1:2})$ for some $a > 0$ and arrive at
\begin{equation*}
\left( \frac{\delta \sigma}{a \abs{\gamma^T\beta_{1:2}}} + 1 \right) a (\gamma^T\beta_{1:2}) = \gamma^T\beta_{1:2}  \,\,\,\,\,\,\,\, \implies \,\,\,\,\,\,\, a = 1 - \frac{\delta \sigma}{\abs{\gamma^T\beta_{1:2}}}. 
\end{equation*}
If this formula yields a negative result, this is a reflection of the fact that the true minimizer $\tilde p_{1:2}$ is orthogonal to $\gamma$, whereupon this derivation is invalid, since the function being minimized is not differentiable at the minimizer, but this can be accounted for by simply setting $a = 0$. Thus we have
\begin{equation*}
a = \max \left\{ 1 - \frac{ \delta \sigma}{\abs{\gamma^T\beta_{1:2}}}, 0 \right\}.  
\end{equation*}
Finally, plugging this back into \eqref{eq:pDerivative2} gives the final update rule
\begin{equation} \label{eq:pUpdate1}
(p_{1:2})_j^{k+1} = \left[\max\left\{0,1 - \frac{\delta \sigma}{\abs{\gamma^T\beta_{1:2}}}\right\} - 1\right](\gamma^T\beta_{1:2})\gamma + \beta_{1:2}.
\end{equation}
One can resolve the minimization for $p_3$ in a similar manner and arrive at
\begin{equation} \label{eq:pUpdate2}
(p_3)_{j}^{k+1} = \max \left\{0, 1 - \frac{\delta\sigma  W}{\abs{\beta_3}}\right\} \beta_3.
\end{equation}

Next we resolve the $k+1$ iterate for the path $x$. First, we solve \eqref{eq:x0mincar} with the initial data function $g(x) = \frac 1 2\abs{x-x_f}^2$. Note that these are vector quantities; in a slight abuse of notation we are letting $x = (x_1,x_2,x_3) = (x,y,\theta)$. Then \eqref{eq:x0mincar} can be written $$\RV{x_0^{k+1} = \argmin_{\tilde x}\left\{\frac\tau2 \vert\tilde x - x_f\vert^2 + \frac 1 2 \vert\tilde x - (x_0^{k} + \tau p_{1}^{k+1})\vert^2 \right\}.} $$ \RV{One easily finds that this has minimizer} 
\begin{equation}  
\RV{x_0^{k+1} = \frac{x_0^k + \tau(x_f + p^{k+1}_1)}{1+\tau}.}
\label{eq:x0k+1}
\end{equation}
To update $x^{k+1}_{j}$ for $j = 1,\ldots, N-1$, we note that our Hamiltonian  $H(x,p)$ given in \eqref{eq:carHam} does not depend on $x_{1:2}$. Because of this, it is trivial to resolve the first and second coordinate in the minimization \eqref{eq:xmincar}:  
\begin{equation} \label{eq:carSpaceUpdate}
(x_{1:2})_j^{k+1} = (x_{1:2})_j^{k} - \tau \left((p_{1:2})_j^{k+1} - (p_{1:2})_{j+1}^{k+1}\right)
\end{equation}
By contrast, one cannot resolve the third coordinate $(x_3)_j^{k+1}$ of the minimizer in \eqref{eq:xmincar} analytically, as the solution to the minimization problem is given in terms of a transcendental equation. Thus we approximate using gradient descent. That is, we set $\theta^* = (x_3)_j^{k+1}$ and perform a few iterations of 
\begin{equation} \label{eq:carGrad1}
\theta^* = \theta^* - \eta h'(\theta^*)
\end{equation}
where $h: \mathbb R \to \mathbb R$ is defined
\begin{equation} \label{eq:carGrad2}\begin{split} 
h(\theta) &= -\delta \tau|(p_1)^{k+1}_{j} \cos(\theta) + (p_2)^{k+1}_j \sin(\theta)| \\ &\hspace{1.5cm}+ \frac 1 2 \big( \theta - ((x_3)^k_j - \tau( (p_{3})^{k+1}_j - (p_3)^{k+1}_{j+1})) \big)^2. \end{split}
\end{equation} We then assign $(x_3)^{k+1}_j = \theta^*.$ Here $\eta > 0$ is the gradient descent rate.  For our simulations, we choose $\eta = 0.15$, and performed three steps of gradient descent at each iteration, though the method also worked with other choices.

Finally, we update the $z$ values:
\begin{equation*}
z_j^{k+1} = x_j^{k+1} + \kappa (x_j^{k+1} - x_j^{k}).
\end{equation*}

\subsection{Implementation of algorithm \ref{alg:1} for the Dubins airplane} \label{plane_updates}
For the Dubins airplane, we use the Hamiltonian \eqref{eq:HJBPlane} which can be written
$$H(x,p) = -p_1 \cos(x_3) - p_2 \sin(x_3) + W_z|p_3| + W_{xy}|p_4|.$$ Here our notation is  $x = (x_1,x_2,x_3,x_4) = (x,y,z,\theta)$ and $p = (p_1,p_2,p_3,p_4) = (u_x,u_y,u_z,u_\theta).$ The update rules for this Hamiltonian are virtually unchanged from those of section \ref{car_updates}, since the minimization for $p_{1:2}, p_3, p_4$ can be resolved individually. 

In fact, the minimization for $p_{1:2}$ is simpler, since there are no absolute values on the first two terms in this Hamiltonian. Define $$\gamma = (\cos((x_4)_{j}^k), \sin((x_4)^k_{j})), \,\,\,\,\,\, \beta = p_j^{k} + \sigma \left(z_j^{k} - z_{j-1}^{k} \right).$$ Note that these are analogous to the definitions of $\gamma$ and $\beta$  in section \ref{car_updates}. Then the update rule for $(p_{1:2})^{k+1}_j$ is 
$$(p_{1:2})_j^{k+1} = \beta_{1:2} +  \delta\sigma \gamma.$$
The updates for the components $p_3$ and $p_4$ have the same structure as before, but with their respective constraints on the angular velocities, $W_z$ for $z$ and $W_{xy}$ for $\theta$. That is, 
$$(p_3)_{j}^{k+1} = \max \left\{0, 1 - \frac{\delta \sigma W_z}{|\beta_3|}\right\} \beta_3,$$
$$(p_4)_{j}^{k+1} = \max \left\{0, 1 - \frac{\delta \sigma W_{xy}}{|\beta_4|}\right\} \beta_4.$$

The updates for the $x$ vector are also similar to those in section \ref{car_updates}.  The update for $x^{k+1}_0$ remains formally the same, though all vectors are four dimensional. The update for the spatial coordinates $x_{1:3}$ is modified only in that in includes the $x_3 = z$ coordinate: 
\begin{equation} \label{eq:planeSpaceUpdate}(x_{1:3})_j^{k+1} = (x_{1:3})_j^{k} - \tau \left((p_{1:3})_j^{k+1} - (p_{1:3})_{j+1}^{k+1}\right).
\end{equation}
The update for $(x_4)^{k+1}_j$ is slightly different (again: simpler). We perform some fixed number of steps of the gradient descent 
\begin{equation} \label{eq:planeGrad1}\theta^*=\theta^* - \eta h(\theta^*) \end{equation} where in this case $h: \R \to \R$ is defined \begin{equation} \label{eq:planeGrad2}\begin{split} h(\theta) &= \delta \tau \big( (p_1)^{k+1}_j \cos(\theta) + (p_2)^{k+1}_j\sin(\theta) \big) \\ &\hspace{1.5cm}+\frac 1 2 \big( \theta - ((x_4)^k_j - \tau( (p_{4})^{k+1}_j - (p_4)^{k+1}_{j+1})) \big)^2 
\end{split} \end{equation} and then assign $(x^4)_{j}^{k+1} = \theta^*.$ The $z$ update is the same as in section \ref{car_updates}.

\subsection{Implementation of algorithm \ref{alg:1} for the Dubins submarine} \label{submarine_updates}
The Hamiltonian for the Dubins submarine in \eqref{eq:HJBSub} can be written $$H(x,p) = \abs{p_1 \cos(x_4)\sin(x_5) + p_1 \sin(x_4)\sin(x_5) + p_3 \cos(x_5)} + W \sqrt{A(x_5)p_{4:5}}$$ where the matrix $A$ is given by $$A(x_5) = \begin{bmatrix} \frac{1}{\sin^2(x_5)} & 0 \\ 0 & 1 \end{bmatrix}.$$ Here our state vector is $x = (x_1,x_2,x_3,x_4,x_5) = (x,y,z,\theta,\varphi)$. 

For ease of notation, when resolving \eqref{eq:pmincar}-\eqref{eq:xmincar} for this Hamiltonian at a time step $j$ on iteration $k+1$, we define 
\begin{align*} 
\gamma &= (\cos((x_4)^{k}_j)\sin((x_5)^{k}_j), \sin((x_4)^{k}_j)\sin((x_5)^{k}_j), \cos((x_5)^k_{j}),\\
\beta &= p^k_j + \sigma(z^k_j - z^{k}_{j-1}),\\
\nu &= x^{k}_j - \tau(p^{k+1}_j - p^{k+1}_{j+1}).
\end{align*}

Again, we can resolve \eqref{eq:pmincar} mostly analytically, by noting that the dependence of $H(x,p)$ on $p_{1:3}$ is decoupled from the dependence on $p_{4:5}$. The update for $p_{1:3}$ is very similarly to the update for $p_{1:2}$ in section \ref{car_updates}. We find the rule
$$
(p_{1:3})^{k+1}_j = M(\gamma^T \beta_{1:3})\gamma  + (\beta_{1:3})
$$
where the constant $M$ is given by 
\begin{equation*} 
M = \max \left(0, 1 - \frac{\delta \sigma}{\abs{\gamma^T\beta}}\right) - 1.
\end{equation*} 

The updates for $p_{4:5}$ are a bit more complicated. Specifically, we need to resolve $$(p_{4:5})_j^{k+1} = \argmin_{\tilde p} \{\delta \sigma \abs{A((x_{5})^{k}_j)\tilde p} + \frac 1 2 \vert \tilde p - \beta_{4:5}\vert^2 \}.$$ This is an evaluation of the proximal operator of the function $f(\tilde p) = \abs{A\tilde p}.$  Using methods similar to those in section \ref{car_updates}, one can resolve this almost analytically. The solution is 
\begin{equation} \label{eq:bisect1}(p_{4:5})^{k+1}_j = \begin{cases} 0, & a^* \le 0, \\ \left(I + \frac{ \delta\sigma  W}{a^*} A((x_5)^k_j)^2\right)^{-1} \beta_{4:5}, & \alpha^* > 0\end{cases}\end{equation}
where $I$ is the $2\times 2$ identity matrix and $\alpha^*$ is the unique root of the function
\begin{equation} \label{eq:bisect2}
g(a) = \frac{\frac 1{\sin^4((x_5)^k_j)} \beta_4^2}{\left(\frac {\delta \sigma  W}{\sin^4((x_5)^k_j)} + a\right)^2} + \frac{\beta_5^2}{( \delta \sigma W + a)^2} - 1.
\end{equation} This function is decreasing, approaches $+\infty$ as $a \searrow -\delta \sigma W$, and approaches $-1$ as $a \nearrow +\infty$, so it has a unique root $a^* \in [-\sigma\delta W, \infty).$ We find this root $a^*$ using a bisection method. First, we check the values $g(-\delta\sigma W + 100\ell)$, for $\ell=1,2,3,\ldots$ to determine the interval of length 100 which contains the root, and then apply the standard bisection method to resolve the root to within $10^{-8}$. 

For the $x$ vector, our Hamiltonian is again independent of $x_{1:3}$, while $x_{4:5}$ will need to be approximated with gradient descent. Specifically, 
\begin{equation} \label{eq:subSpaceUpdate}
(x_{1:3})^{k+1}_j = \nu_{1:3},
\end{equation}
and for $x_{4:5}$, we perform a few iterations of 
\begin{equation} \label{eq:subGrad1}(\theta^*,\varphi^*) = (\theta^*,\varphi^*) - \eta \nabla h(\theta^*,\varphi^*) \end{equation}
where \begin{equation} \label{eq:subGrad2}\begin{split} h(\theta,\varphi) &= -\delta \tau \Bigg(\abs{(p_{1})^{k+1}_j \cos(\theta)\sin(\varphi) + (p_{2})^{k+1}_j \sin(\theta) \sin(\varphi) + (p_3)^{k+1}_j \cos(\varphi)} + \\ &\hspace{1.5cm} W \sqrt{\frac{(p_4)^{k+1}_j}{\sin^2(\varphi)} + (p_5)^{k+1}_j}\Bigg) + \frac 1 2((\theta - \nu_4 )^2 + (\varphi - \nu_5)^2), \end{split} \end{equation} and then assign $(x_{4:5})^{k+1}_j = (\theta^*,\varphi^*).$ 

As mentioned before, anywhere one sees a $\sin^2(\varphi)$ in a denominator, we replace it with $\sin^2(\varphi) + \eps$ where $\eps = 10^{-10}$, to avoid division by zero. Once again, the $z$ updates are the same. 

\subsection{Adding Obstacles Computationally}
One final matter to address is how to computationally account for obstacles. As mentioned in section \ref{sec:TimeHorizonAndObstacles}, we can account for obstacles without using boundary conditions by modifying the Hamilton-Jacobi equation by multiplying the Hamiltonian by the indicator function of the free space. Thus we are actually solving \begin{equation} \label{eq:HJWithObs}
u_t + O(x,t)H(x,\nabla u) = 0,  \,\,\,\,\,\, u(x,0) = g(x),
\end{equation} where $O(x,t) = 0$ when $x \in \Omega_{\text{obs}}(t)$ and $O(x,t) = 1$ otherwise. 

However, this raises the question of how to efficiently determine when $x \in \Omega_{\text{obs}}(t)$ for a given point $x$ and time $t.$ To make this possible, we restrict ourselves to the case of obstacles which are disjoint collections of balls. Thus at any time $t > 0$, we have $$\Omega_{\text{obs}}(t) = \bigcup^{L}_{\ell=1} B(x_\ell(t),r_\ell(t)).$$ In this case, it is computationally efficient to check if $x$ is in a obstacle at time $t$ by checking its distance to each of the centers $x_\ell(t).$ [Note, there is another small abuse of notation here: the centers of the obstacles will be \emph{spatial} coordinates only, whereas our state vector has spatial and angular coordinates.]

Of course, in application, it will almost never be the case that the set of obstacles is a finite collection of balls. In this case, we run the following greedy algorithm to iteratively fill the obstacle with disjoint balls of the largest possible radius.
\begin{itemize} \item[(0)] Given a set of obstacles $\Omega_0 = \Omega_{\text{obs}}(t)$ and a minimal radius $r_\text{min}$, set $\ell = 1$ and do the following
    \item[(1)] Compute the signed distance function $d_{\ell-1}(x)$ to $\partial \Omega_{\ell-1}$ using the level set method \cite{OsherFedkiw}, fast marching method \cite{SethianLevel}, parallel redistancing method \cite{Redist1,Redist2}, or another method of your choice. We choose the distance which is \emph{positive} inside $\Omega_{\ell-1}$ and negative outside the set.
    \item[(2)] Set $x_{\ell} = \argmax_x d_{\ell-1}(x)$, $r_\ell = d_{\ell-1}(x_\ell),$ and $B_\ell = B(x_\ell,r_\ell).$
    \item[(3)] Set $\Omega_\ell = \Omega_{\ell-1} \setminus B_\ell$.
    \item[(4)] If $r_\ell < r_{\text{min}}$,  break the loop. Otherwise, increment $\ell$ and return to (1). 
\end{itemize}

\begin{figure}[t!]
    \centering
    \includegraphics[width=0.31\textwidth]{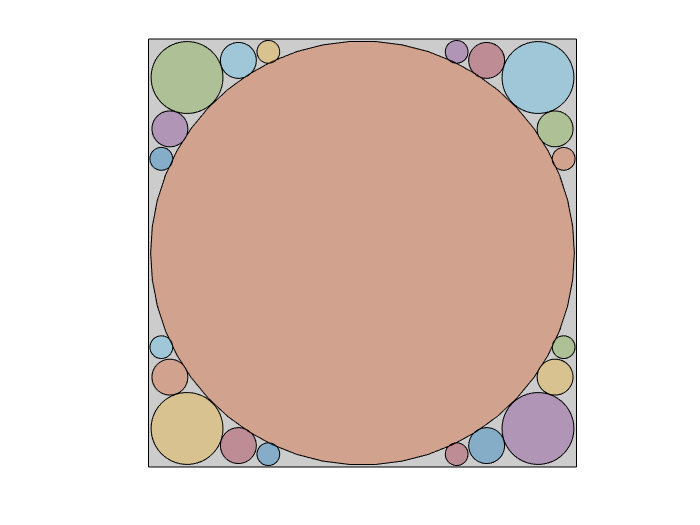}\,
    \includegraphics[width=0.31\textwidth]{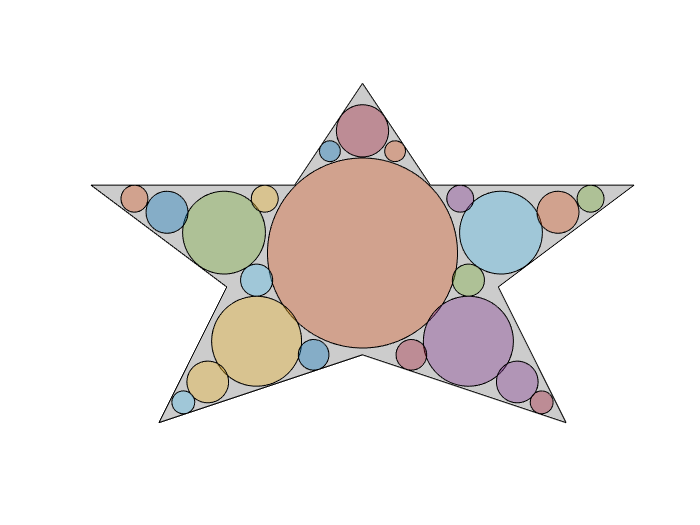}\,
    \includegraphics[width=0.31\textwidth]{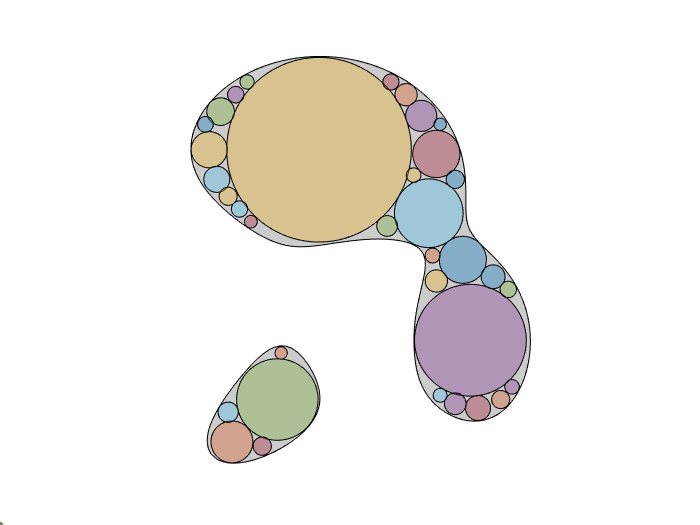}
    \caption{Three obstacles approximated by disjoint circles.}
    \label{fig:obs}
\end{figure}

Doing this results in a collection of disjoint balls which approximate the shape $\Omega_{\text{obs}}(t)$ and have radii larger than the prespecified minimum radius $r_{\text{min}}.$ This is demonstrated in figure \ref{fig:obs}. \RV{The $r_{\text{min}}$ we use is chosen on an \emph{ad hoc} basis. In cases where the obstacles have "long, thin" parts, one may need to choose $r_{\text{min}}$ very small, while for smoother, more regular obstacles, a larger choice of $r_{\text{min}}$ may suffice. In general, covering of surfaces by balls is a problem of interest in computational geometry \cite{balls1,balls2}, and theory regardig recursive subdivision algorithms for circle packing is an active area of research \cite{balls3}. Our algorithm is admittedly quite simple, but is sufficient for our purposes.}

Note that as described, this algorithm is somewhat computationally burdensome, but it is done in preprocessing. As an alternative, if one wishes to do this in real-time (and possibly do away with the need to approximate the obstacle by circles), one could one may be able to embed the algorthms from \cite{Redist1,Redist2} into algorithm \ref{alg:1}, since they also compute the distance function based on Hopf-Lax type formulas. For our purposes, we use the algorithm described above and approximate obstacles with circles. 

As one final note, we point out how adding the obstacle function $O(x,t)$ to the Hamilton-Jacobi equation as in \eqref{eq:HJWithObs} affects the update rules in each of the cases above. When replacing the Hamiltonian $H(x,p)$ with $O(x_s,t)H(x,p)$, where $x_s$ denotes the spatial variables in $x$, the update rules for any of the costate variables do not change in any significant way: the indicator function of the obstacles is simply brought along as a multiplier and appears anywhere where $\delta\sigma$ appears. For example, the update rules for the car given by \eqref{eq:pUpdate1} and \eqref{eq:pUpdate2} become \begin{align*}(p_{1:2})_j^{k+1} &= \left[\max\left\{0,1 - \frac{O((x_{1:2})^k_j,t_j)\delta \sigma}{\abs{\gamma^T\beta_{1:2}}}\right\} - 1\right](\gamma^T\beta_{1:2})\gamma + \beta_{1:2}, \\
(p_3)_{j}^{k+1} &= \max \left\{0, 1 - \frac{O((x_{1:2})^k_j,t_j)\delta\sigma  W}{\abs{\beta_3}}\right\} \beta_3, \end{align*} which are the exact same as before aside from the inclusion of $O(x_s,t)$ evaluated at the current spatial coordinates $(x_{1:2})^k_{j}$ and time step $t_j$. The updates for the airplane or submarine are modified in analogous ways. 

Likewise, because the obstacle function $O(x_s,t)$ depends only on the spatial variables in the state vector, the updates for the angular variables (described by equations \eqref{eq:carGrad1}, \eqref{eq:carGrad2} for the car, by equations \eqref{eq:planeGrad1}, \eqref{eq:planeGrad2} for the airplane, and by equations \eqref{eq:subGrad1}, \eqref{eq:subGrad2} for the submarine) remain virtually unchanged, except once again for the inclusion of the obstacle function as a multiplier wherever $\delta \tau$ appears. For example, for the car, we simply need to change the function $h$ in \eqref{eq:carGrad2} to \begin{align*} 
h(\theta) &= -\delta \tau O((x_{1:2})^{k}_j,t_j)|(p_1)^{k+1}_{j} \cos(\theta) + (p_2)^{k+1}_j \sin(\theta)| \\ &\hspace{1.5cm}+ \frac 1 2 \big( \theta - ((x_3)^k_j - \tau( (p_{3})^{k+1}_j - (p_3)^{k+1}_{j+1})) \big)^2, \end{align*} which is the exact same as before aside from the inclusion of $O(x_s,t)$. Analogous modifications are made for the plane and submarine. 

The significant change comes in the updates of the spatial components which are given for the car, airplane, and submarine by equations \eqref{eq:carSpaceUpdate},\eqref{eq:planeSpaceUpdate}, and \eqref{eq:subSpaceUpdate} respectively. The minimization problems for these updates can no longer be resolved analytically, so they approximate the solutions using gradient descent, as was done with the angular variables before. In each case, we need to resolve $$(x_s)^{k+1}_j = \argmin_{\overline x} \left\{-\delta \tau O(\overline x_s,t)H(x_a,p) + \frac 1 2 \|\overline x - \nu\|^2\right\}$$ where $\nu = (x_s)^k_{j} - \tau((p_s)^{k+1}_j - (p_s)^{k+1}_{j+1})$. To reiterate, we are using $x_s$ to denote the spatial variables of $x$ (and $p_s$ to denote the corresponding co-state variables), and we introduce $x_a$ to denote the angular variables. Note that in our three models above, the Hamiltonian depends only on the angular variables $x_a$.  To approximate the solution of the minimization problem, we use the gradient descent $$x_s^* = x_s^* - \eta \nabla h(x_s^*)$$ where $h$ is defined $$h(x_s) = -\delta \tau O( x_s,t_j)H((x_a)^k_j,p^{k+1}_j) + \frac 1 2 \vert x_s - \nu\vert^2$$ and then set $(x_s)^{k+1}_j = x_s^*.$ In doing this, we will need to evaluate $\nabla O(\cdot,t)$. Recall, this function is the indicator function of the free space, which is discontinuous across the boundary of the obstacles. Accordingly, we approximate the indicator function of the free space by   
 $$O(x_s,t) \approx \frac{1}{2} + \frac{1}{2}\tanh(-100d(x_s,t))$$ where $d(x_s,t)$ is the signed distance function to the boundary of the obstacles at time $t$ (positive inside the obstacles). This definition gives a smooth function which is approximately zero inside the obstacles and approximately 1 in the free space. Here the $100$ in the definition is an arbitrary large number. We can quickly evaluate the distance function when we are approximating the obstacles by disjoint circles. This also gives an efficient method for computing the gradient of the distance function when $x_s$ is outside the obstacles: it is the unit vector pointing from $x_s$ to the center $x_\ell(t)$ of nearest ball to $x_s$ in the collection of balls which approximates the obstacles. 

 Note that the updates for $x^{k+1}_0$---which are given by equation \eqref{eq:x0k+1} for all three models---do not change because the Hamiltonian does not appear in these updates. 

\section{Simulations \& Discussion} \label{results}

In this section, we present the results of several simulations which demonstrate the efficacy of our models. In all simulations, we choose somewhat arbitrary and synthetic data, which we list for the individual simulations. The simulations were performed on a laptop computer with an AMD Ryzen 7 7735U processor running at a maximum of 4.75 GHz, and 16GB of RAM. Besides anything described above, no further efforts were made to optimize the algorithms, though certain parts are parallelizable, and the resolution of the minimizers used to update the state variables could likely be done more efficiently. Even so, our algorithms are quite efficient as reported below. 

In all cases we fix the time step $\delta = 0.1$ and discretize the paths into $N = t/\delta$ steps, where $t$ is the time chosen for the simulation as seen in algorithm \ref{alg:1}. Whenever it is necessary to use gradient descent to perform the updates described in sections \ref{car_updates}-\ref{submarine_updates}, we perform 3 steps of gradient descent with a descent rate of $\eta = 0.15$. We also fix $\sigma = 0.5, \tau = 0.5$ and $\kappa = 1$. In each case, we continue the iteration until the maximum change in any coordinate of $x$ or $p$ is less than $\text{TOL} = 10^{-3}$. This tolerance may appear fairly large, but there was no appreciable difference in the resolved optimal trajectories when this tolerance was decreased to $10^{-8}$. \RV{In each case, we use a maximal iteration count of $K = 100000$ but this maximal count was never reached when using the baseline parameter values. When we report clock times and iteration counts for the simulations below, these are averaged over 50 trials. Average clock time is then rounded to the nearest tenth of a second and average iteration count is rounded to the nearest integer. }

\begin{figure}[t!]
    \centering
    \includegraphics[width=0.4\textwidth,trim = 50 10 10 0,clip]{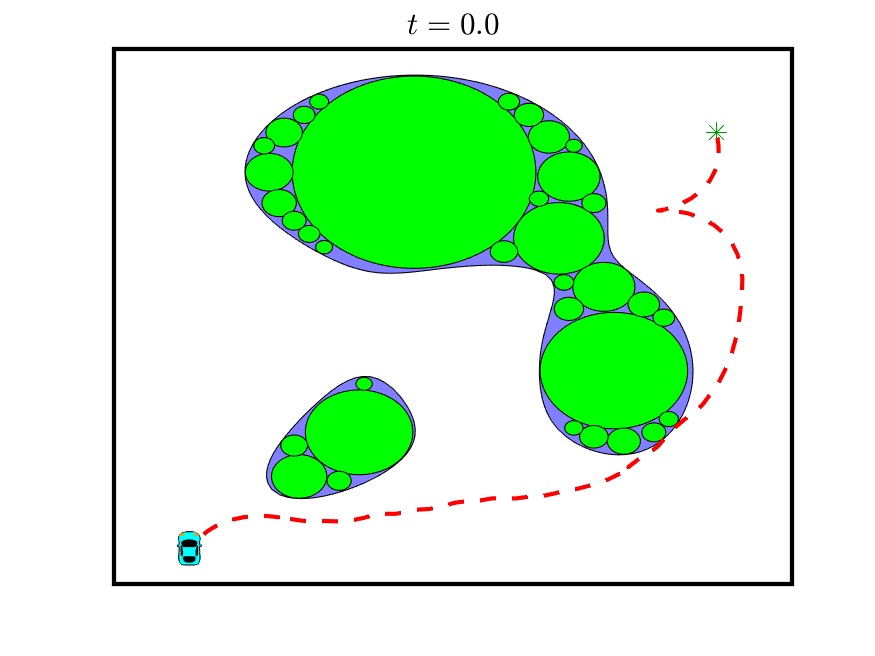}\,
    \includegraphics[width=0.4\textwidth,trim = 50 10 10 0,clip]{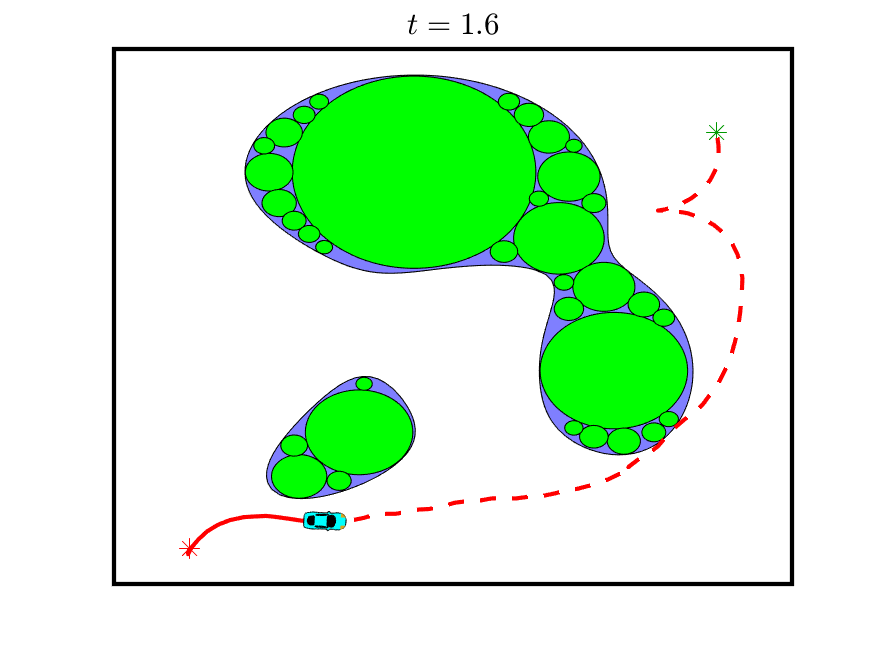}\\
    \includegraphics[width=0.4\textwidth,trim = 50 10 10 0,clip]{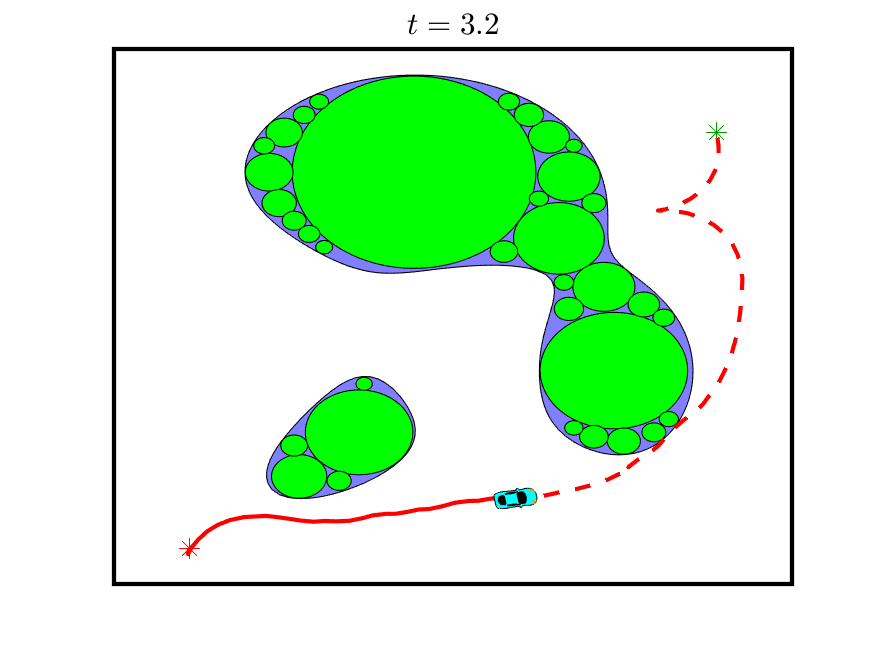}\,
    \includegraphics[width=0.4\textwidth,trim = 50 10 10 0,clip]{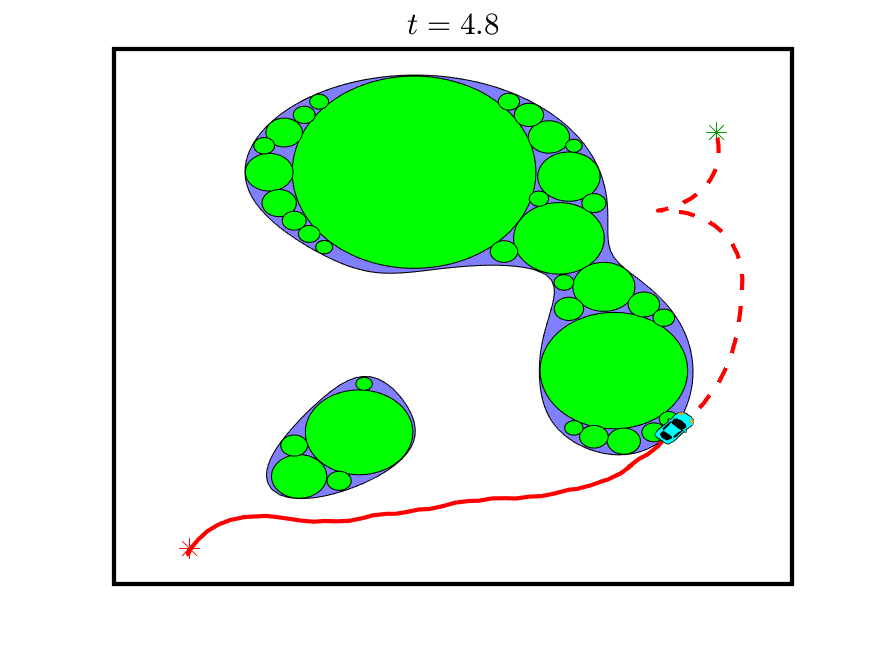}\\
    \includegraphics[width=0.4\textwidth,trim = 50 10 10 0,clip]{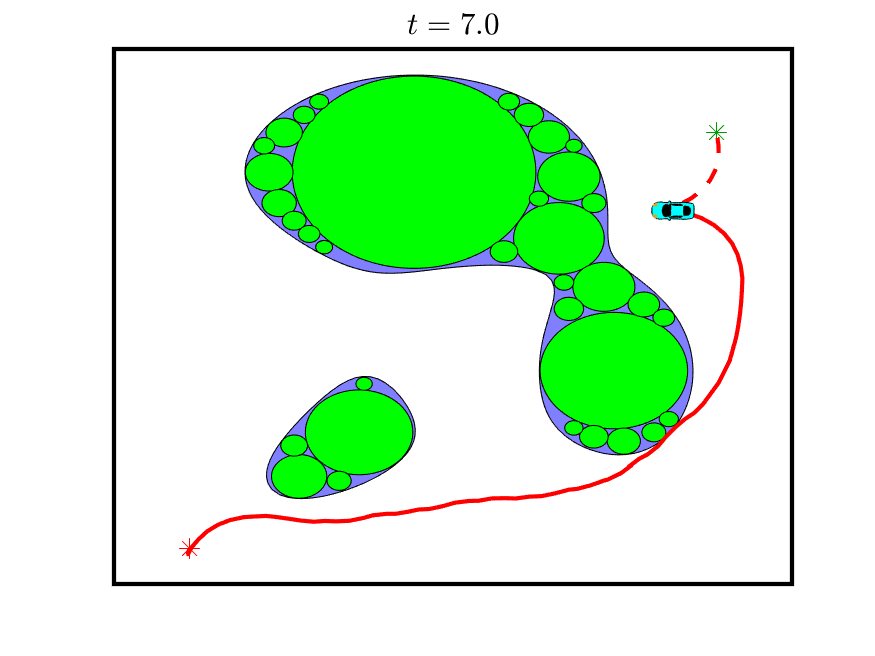}\,
    \includegraphics[width=0.4\textwidth,trim = 50 10 10 0,clip]{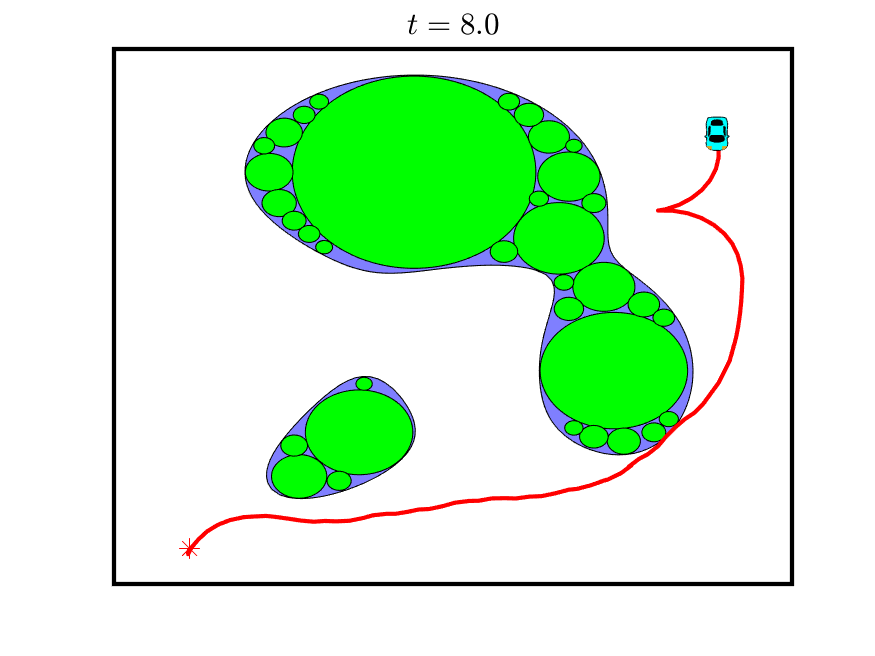}
    \caption{A car navigating around obstacles to end point in the top right of each frame. The obstacles are the blue regions, and the collection of green circles is their computational representation. }
    \label{fig:carS}
\end{figure}

\begin{figure}[t!]
    \centering
    \includegraphics[width=0.4\textwidth,trim = 50 10 10 0,clip]{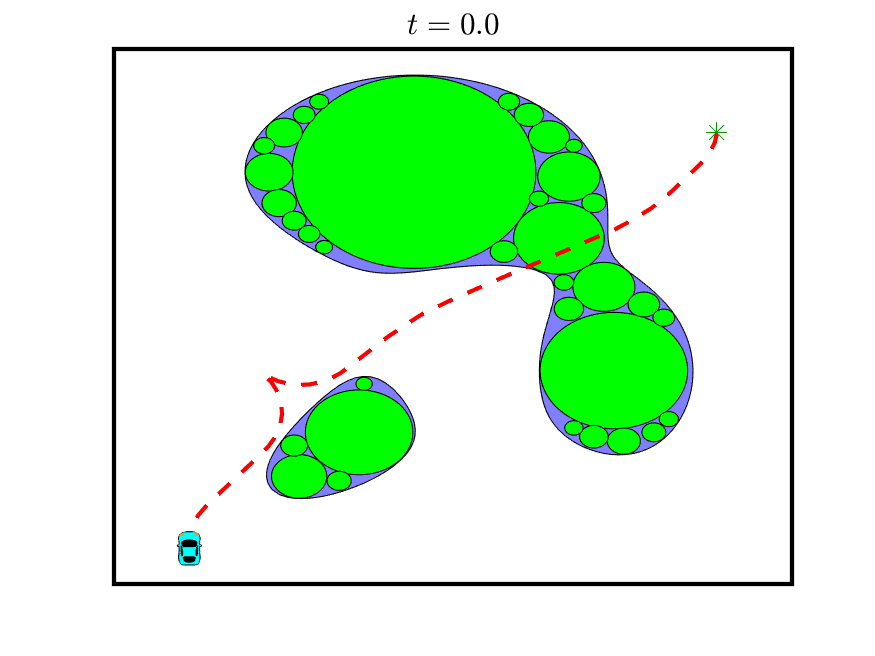}\,
    \includegraphics[width=0.4\textwidth,trim = 50 10 10 0,clip]{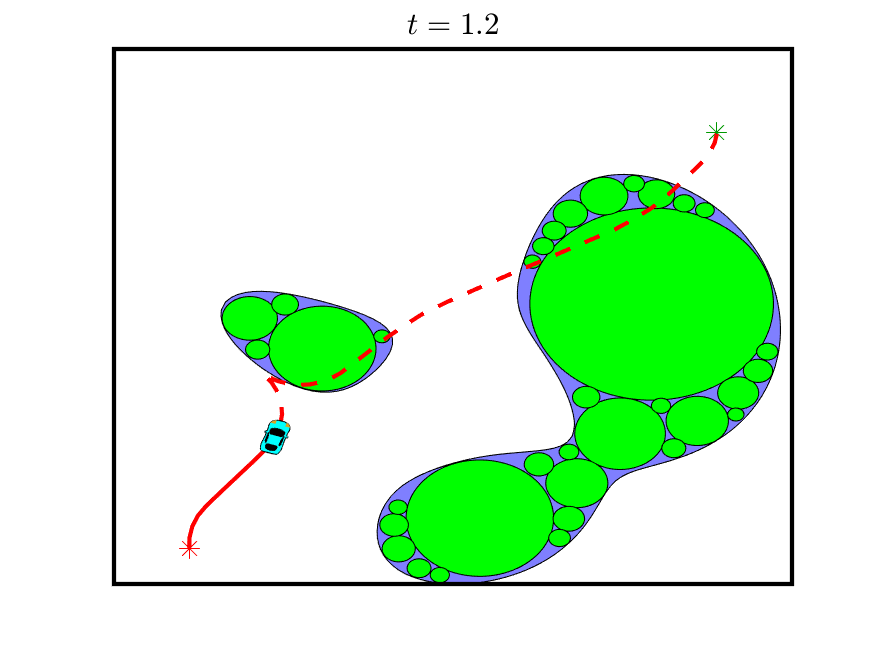}\\
    \includegraphics[width=0.4\textwidth,trim = 50 10 10 0,clip]{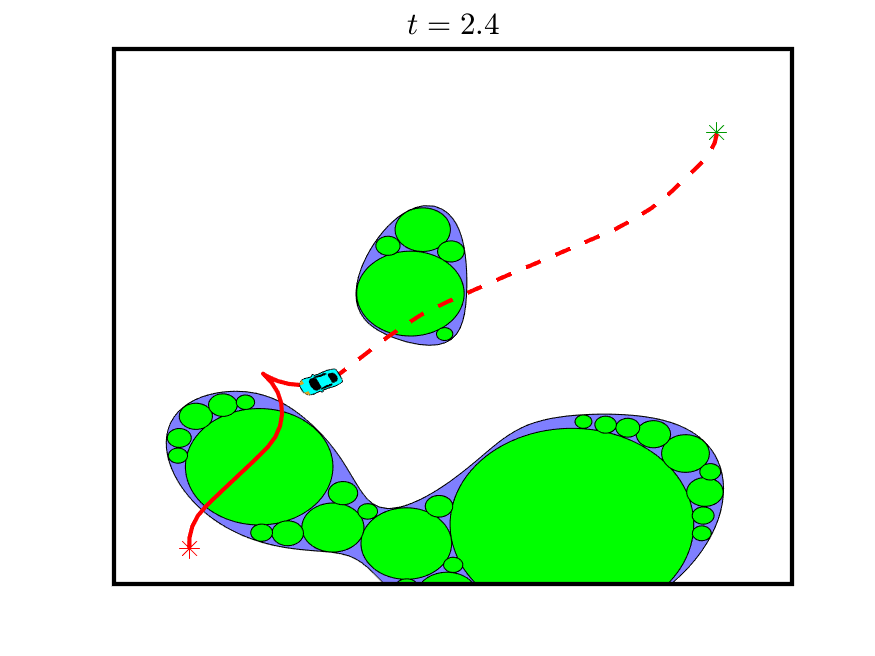}\,
    \includegraphics[width=0.4\textwidth,trim = 50 10 10 0,clip]{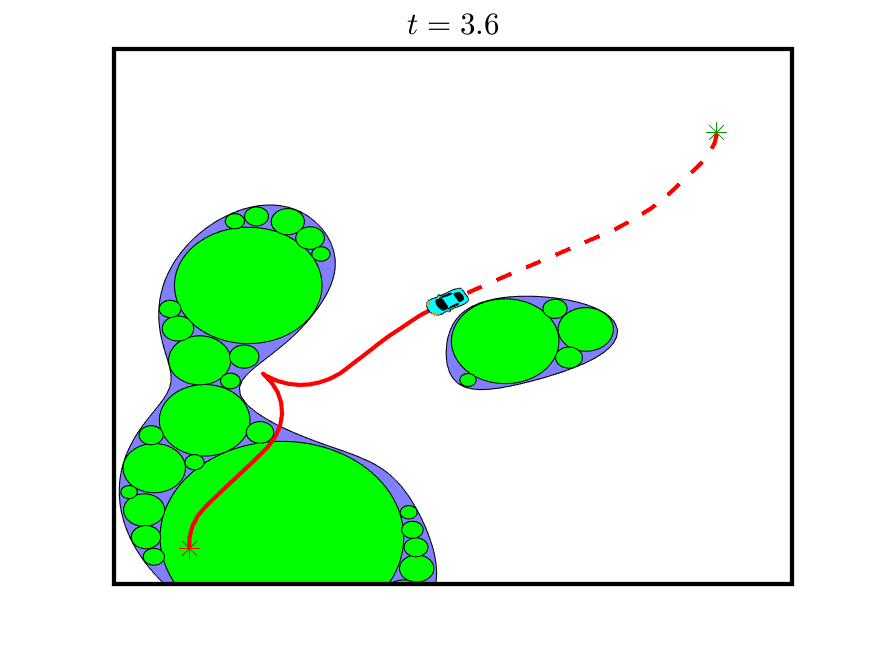}\\
    \includegraphics[width=0.4\textwidth,trim = 50 10 10 0,clip]{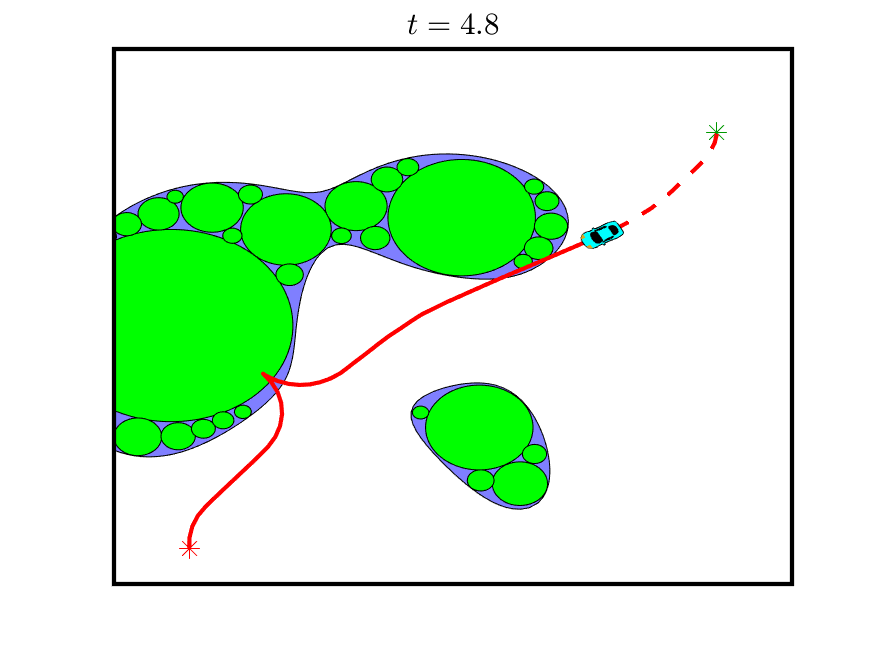}\,
    \includegraphics[width=0.4\textwidth,trim = 50 10 10 0,clip]{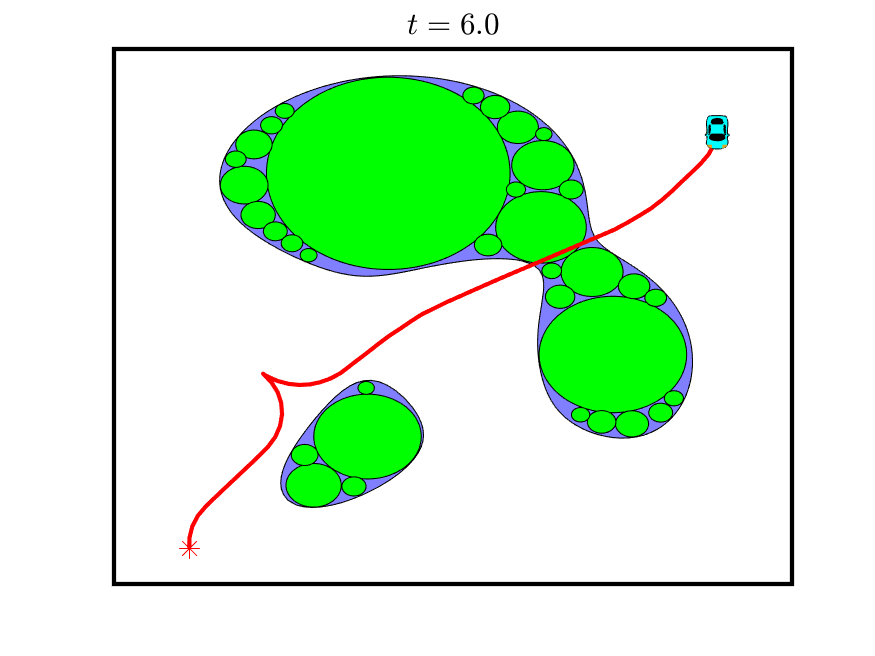}
    \caption{A car navigating around the same obstacles as in figure \ref{fig:carM}, but in this case the obstacles begin rotating clockwise about the center of the domain.}
    \label{fig:carM}
\end{figure}

Our first two simulations, the results of which are shown in figure \ref{fig:carS} and figure \ref{fig:carM}, show a car beginning in the configuration $(x_0,y_0,\theta_0) = (-1.5,1.5,\pi/2)$ and navigating around obstacles to the point $(x_f,y_f,\theta_f) = (2,2,3\pi/2)$ which are in the bottom-left and top-right (respectively) of the domain used to plot the results. In both cases, the maximum angular velocity is $W = 2$. The obstacles in both figures have the same shape, but in figure \ref{fig:carS} the obstacles remain stationary, whereas in figure \ref{fig:carM}, the obstacles rotation clockwise around the origin at a constant rate of 1 radian per unit time. In the former figure, the time required to reach the endpoint is $T = 8$, whereas, when the obstacles rotate out of the way in the latter figure, the car can reach the endpoint by time $T = 6$.  As mentioned in section \ref{sec:TimeHorizonAndObstacles}, the time-horizon $T$ that is chosen for the problem does actually matter here. Theoretically, for any $T$ chosen large enough, there is an optimal path from the initial point to the end point. However, empirically if $T$ is chosen too large the algorithm will require longer to converge. As chosen for each of these simulations, the time horizons of $T = 8$ and $T = 6$, respectively, are very close to the best possible travel time. For the moving obstacles example (figure \ref{fig:carM}, the algorithm was able to resolve optimal paths in an average 1.6 seconds using an average of 1748 iterations. If $T$ is chosen too small, so that no path requiring time less than $T$ to traverse can reach the endpoint, then the algorithm will often fail to converge, or converge to a path which is not meaningful.  In these examples if the tolerance for convergence is decreased to $\text{TOL} = 10^{-8}$, the optimal trajectories do not appreciably change, but the clock time increases to roughly 15 seconds, requiring on the order of 10000-15000 iterations for convergence. 

\begin{figure}[t!]
    \centering
    \includegraphics[height=0.22\textheight,trim= 50 0 50 0,clip]{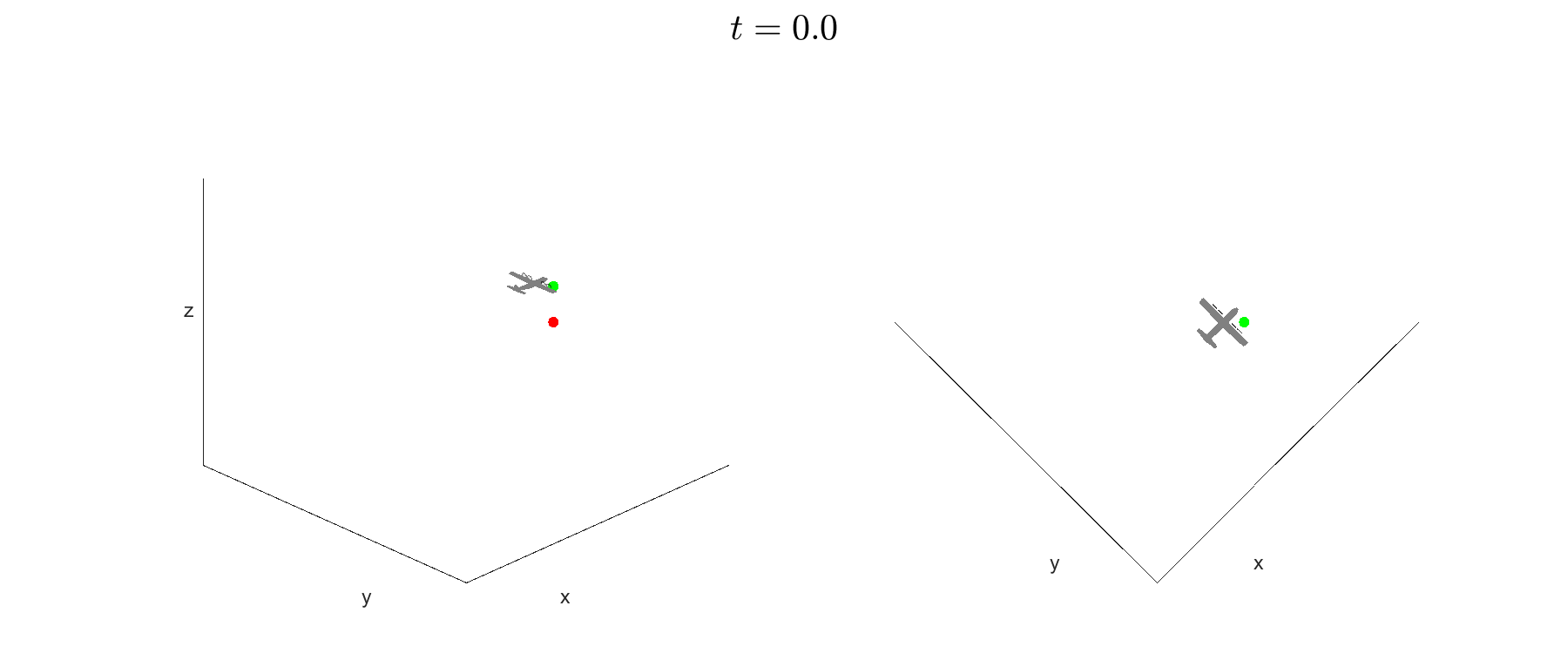}\\
    \includegraphics[height=0.22\textheight,trim= 50 0 50 0,clip]{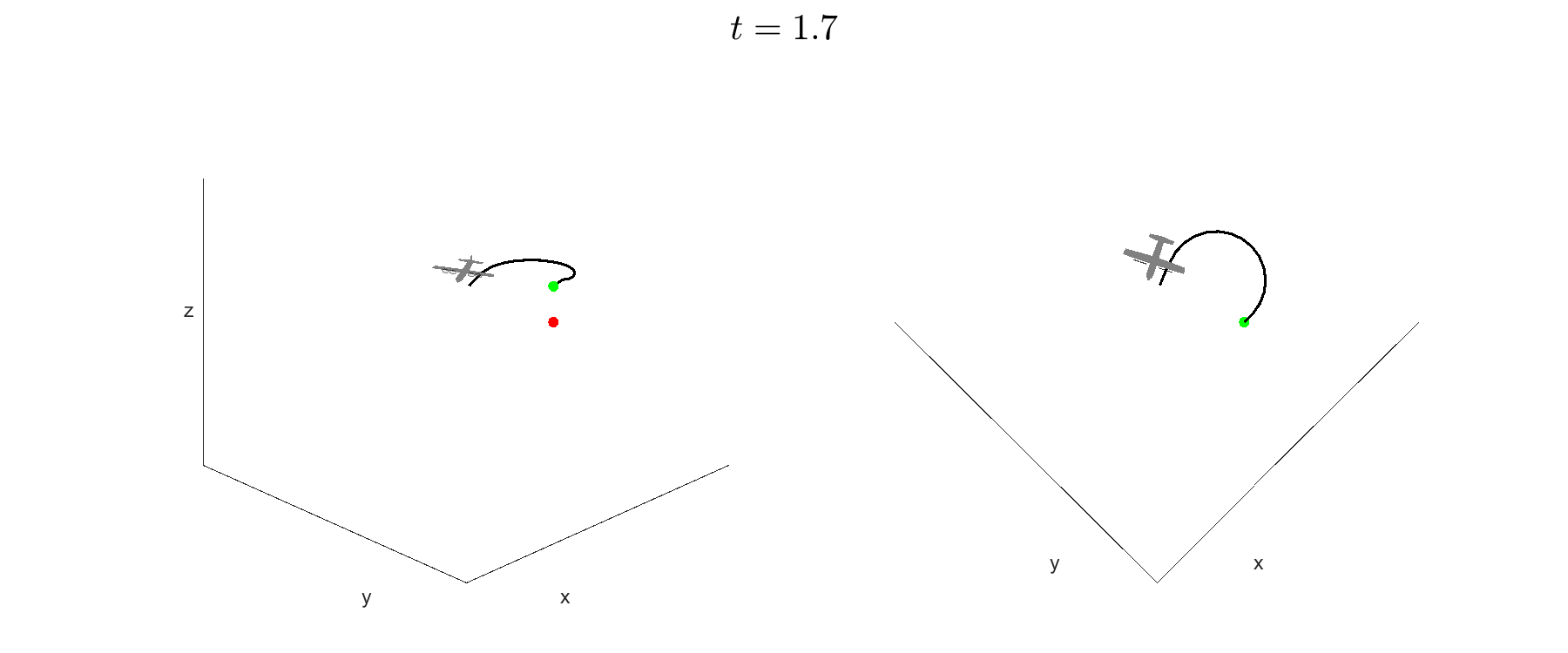}\\
    \includegraphics[height=0.22\textheight,trim= 50 0 50 0,clip]{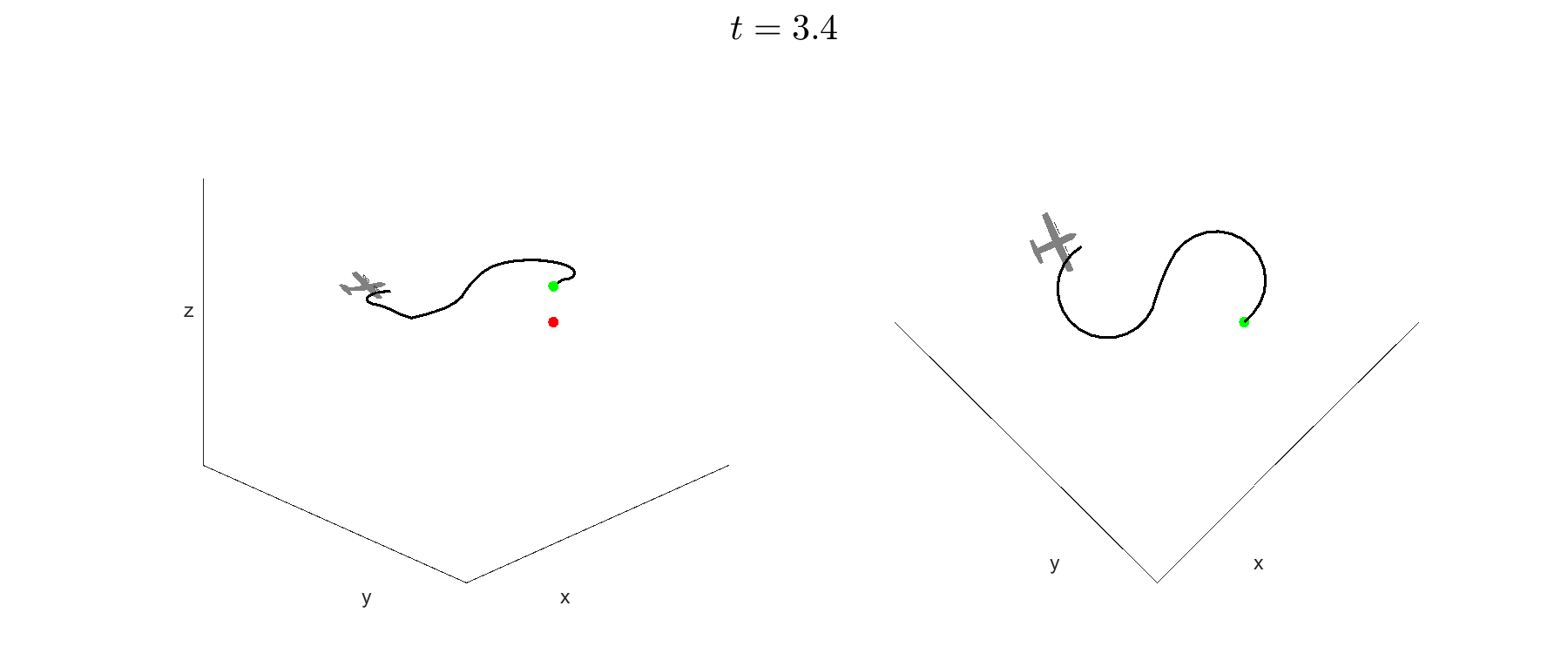}\\
    \includegraphics[height=0.22\textheight,trim= 50 0 50 0,clip]{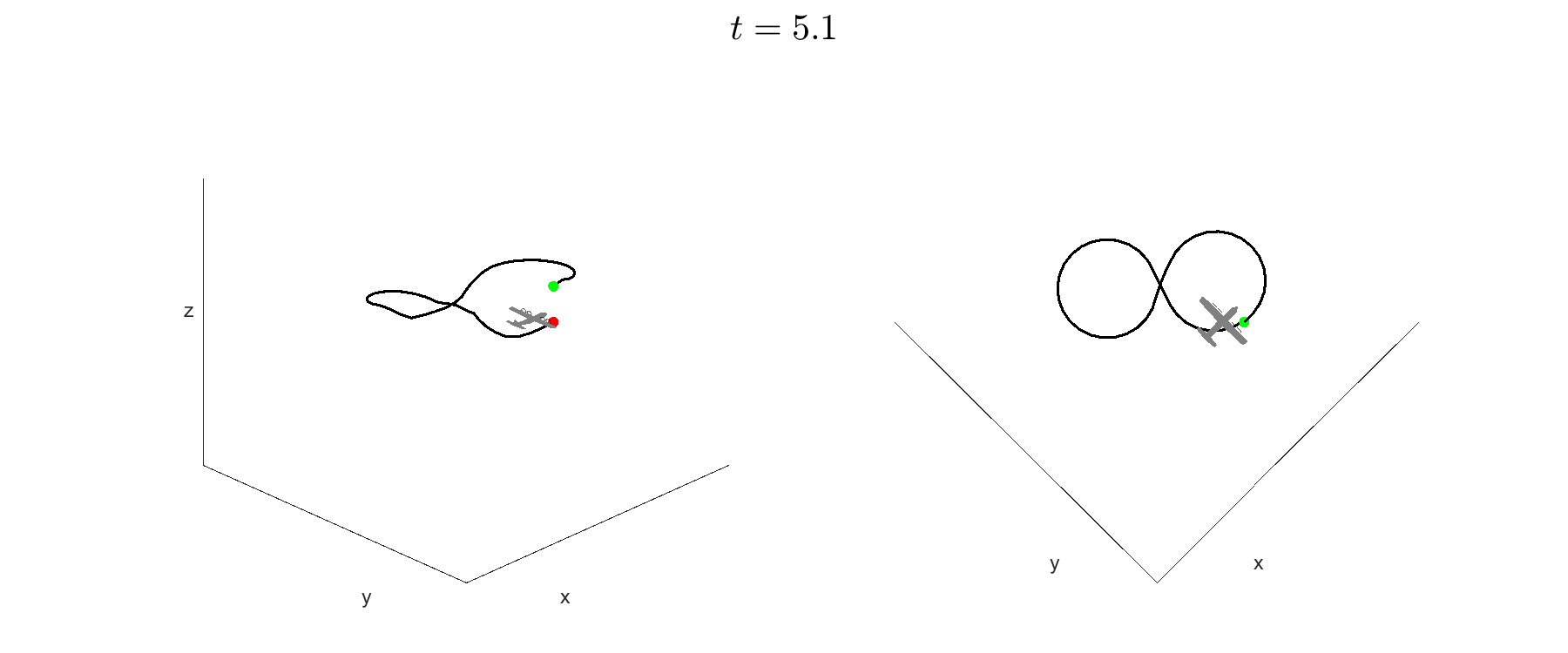}
    \caption{A Dubins airplane circling down for a landing. The 3D view is displayed on the left of each panel, and the projection down to the $xy$-axis is displayed on the right.}
    \label{fig:plane}
\end{figure}

Our next simulation has a Dubins airplane circling down for a landing, which is displayed in figure \ref{fig:plane}. The airplane begins in the configuration $(x_0,y_0,z_0,\theta_0) = (0,0,\tfrac 1 2,0)$ and must end at $(x_f,y_f,z_f,\theta_f) = (0,0,0,0).$ In this case, the maximum angular velocity in the $xy$-plane is $W_{xy} = 2.5$ and the maximum angular velocity in the $z$-direction is $W_z = 0.5$. In order to make its descent, the airplane flies in a perfect figure eight in the $xy$-plane, while descending in the $z$-direction. In figure \ref{fig:plane}, the 3D view is displayed on the left of each panel and the projection down to the $xy$-plane is displayed on the right of each panel. \RV{This simulation required an average of 0.4 seconds of clock time to resolve the optimal path, doing so in an average of 2506 iterations. The updates for the plane are a bit simpler than those for the car (since the plane can only move forward) which explains why more iterations can be performed in less clock time. For this example, there is not much penalty for deceasing the convergence tolerance. Indeed, decreasing the tolerance all the way to $\text{TOL} = 10^{-8}$, the algorithm still usually resolves the optimal path in an average of 1.2 seconds and 7233 iterations (though the ``worst" simulation required 5.4 seconds and 32268 iterations). An interesting note is that the Dubins airplane does not have small time local controllability (STLC), meaning roughly that the plane cannot necessarily reach a point which is with distance $\eps$ of its current configuration in time $O(\eps)$. The value function for a control problem without STLC can be discontinuous, which necessitates special consideration when designing grid based numerical schemes for these problems as discussed in \cite{TakeiTsai2}. Because our method only ever deals with the Hamiltonian and initial data function (and never with approximations to derivatives of the value function), our method is agnostic to whether the dynamics satisfy an STLC condition.}

\begin{figure}[t!]
    \centering
    \includegraphics[height=0.215\textheight,trim= 50 0 50 0,clip]{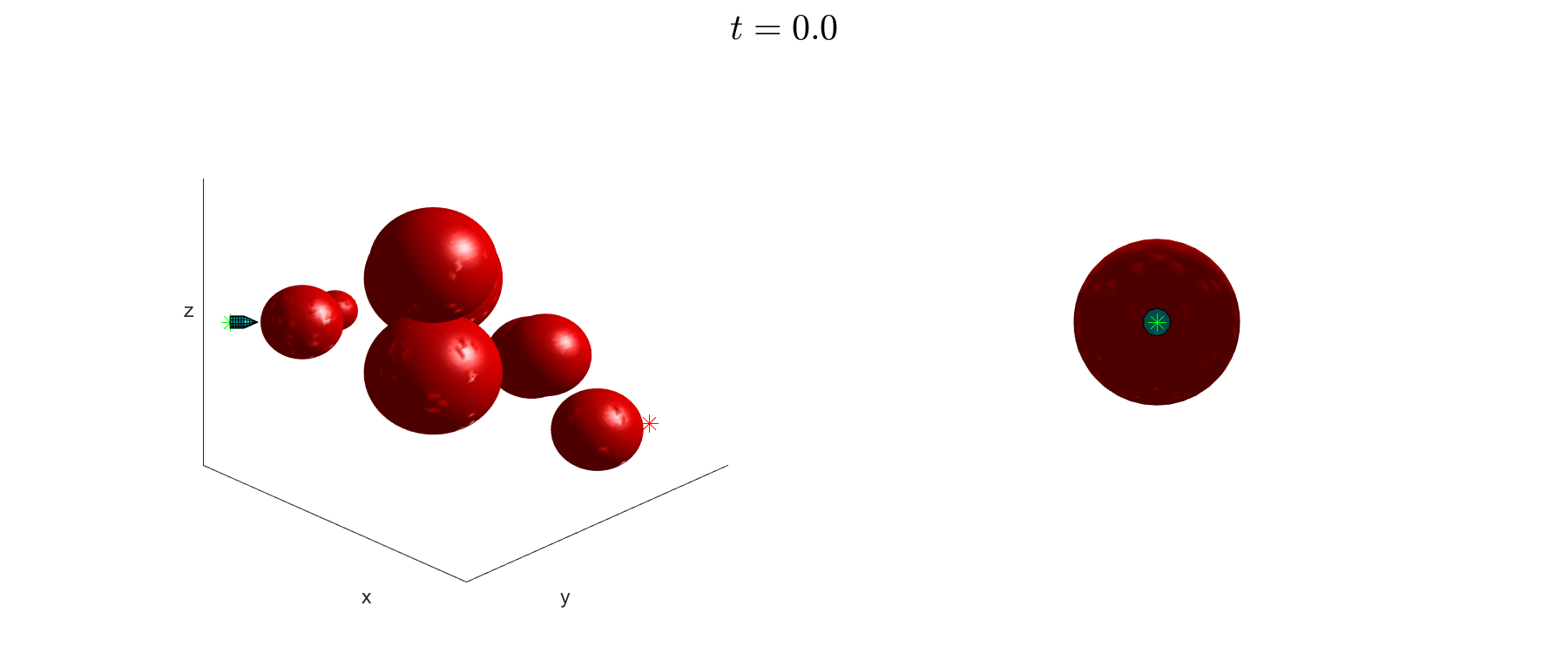}\\
    \includegraphics[height=0.215\textheight,trim= 50 0 50 0,clip]{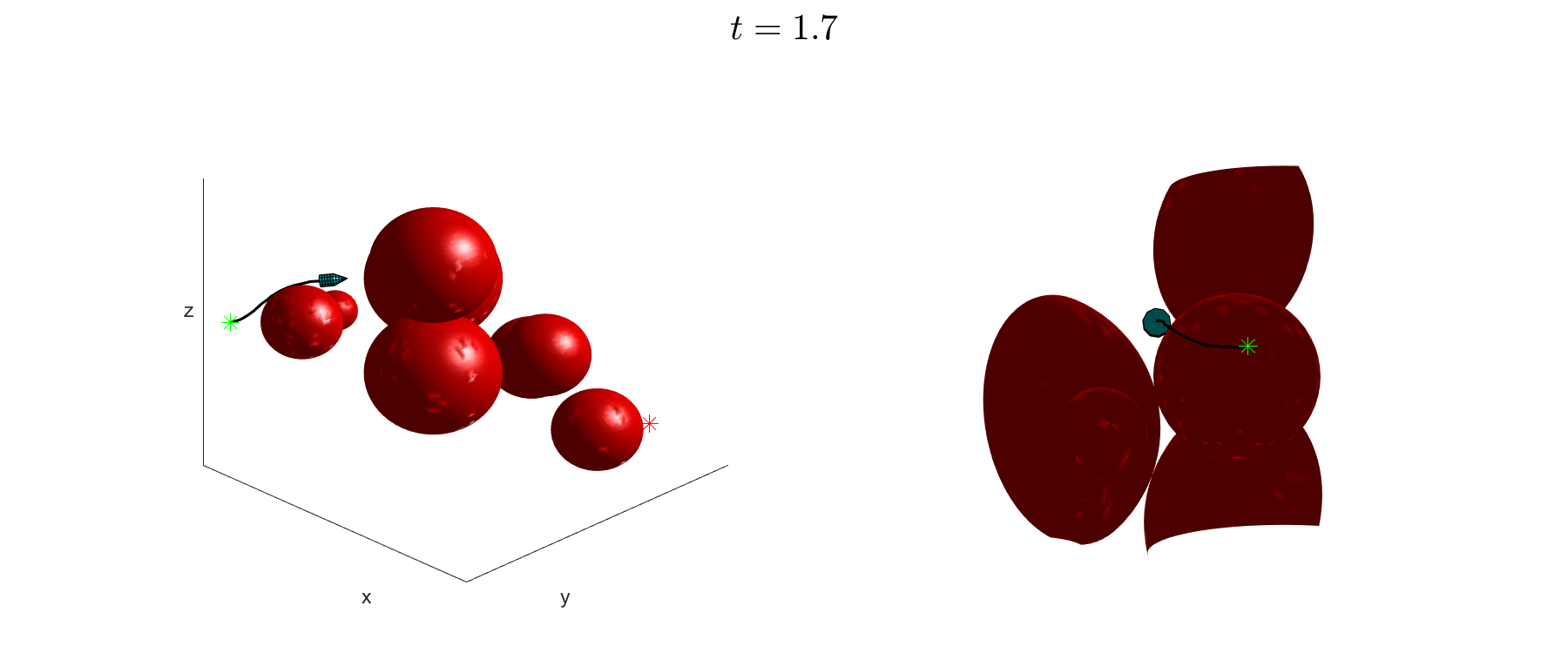}\\
    \includegraphics[height=0.215\textheight,trim= 50 0 50 0,clip]{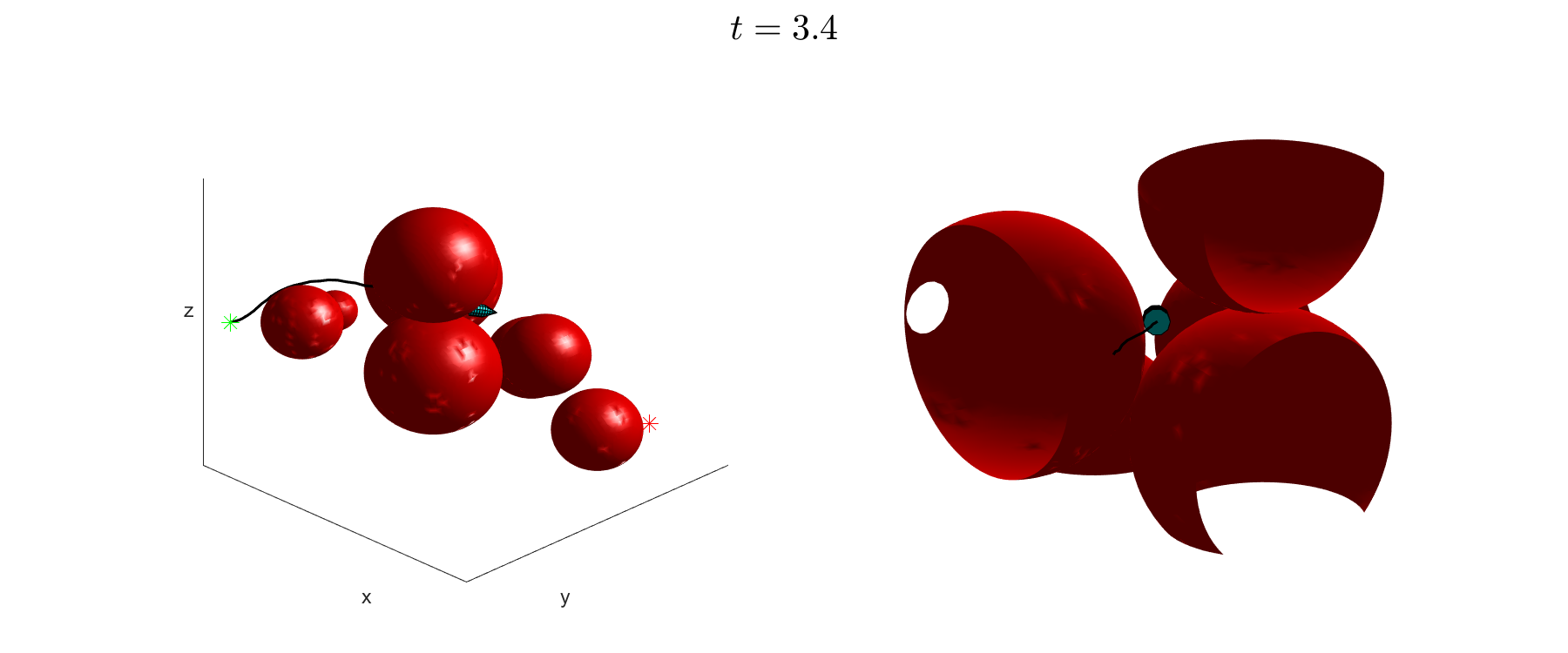}\\
    \includegraphics[height=0.215\textheight,trim= 50 0 50 0,clip]{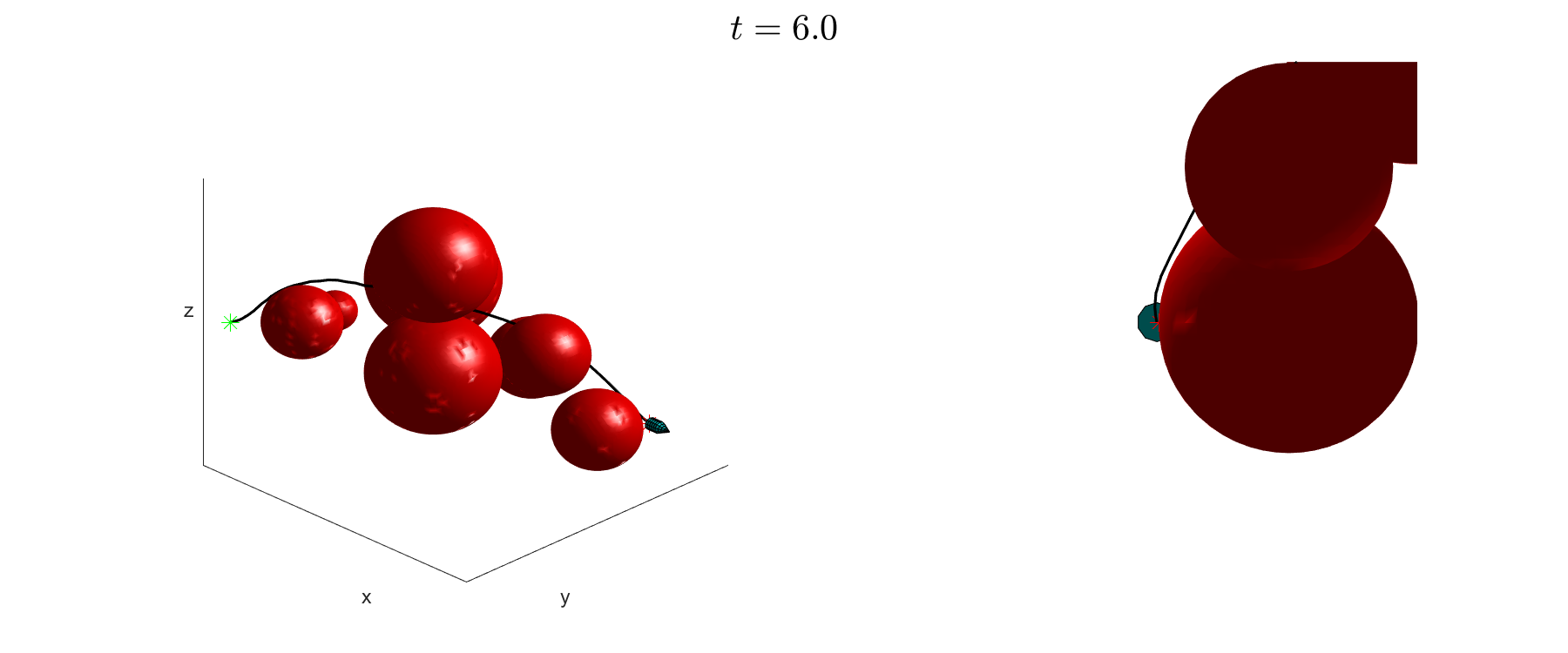}
    \caption{A Dubins submarine navigating through and around obstacles to reach its desired destination. The ``third-person" view is displayed on the left of each panel, while the corresponding ``first-person" view is displayed on the right. The red bubbles are obstacles.}
    \label{fig:sub}
\end{figure}

We next consider a Dubins submarine. Here we take the initial configuration to be $(x_0,y_0,z_0,\theta_0,\varphi_0) = (-1.8,-1.8,0,\pi/4,\pi/2)$ and let the submarine navigate through and around obstacles to $(x_f,y_f,z_f,\theta_f,\varphi_f) = (1.3,1.5,-1.5,0,\pi/2).$ In this case, the maximum angular velocity is $W=2$. The result is displayed in figure \ref{fig:sub}, where the ``third person" view is displayed on the left of each panel and the corresponding ``first person" view is displayed on the right. In figure \ref{fig:sub}, the red bubbles are obstacles. This simulation required an average of 5.0 seconds of clock time to resolve the optimal trajectory in an average of 1936 iterations. Here the iterations are more expensive due to the resolution of $p_{4:5}$ by the bisection method described in equations \eqref{eq:bisect1}, \eqref{eq:bisect2}. Once again, if we decrease the convergence tolerance to $\text{TOL} = 10^{-8}$, the path generated is not appreciably different, but requires roughly 28 seconds of clock time and on the order of 8000 iterations to resolve. 

Finally, we want to empirically study the affect that varying the ambient parameters has on the convergence of the algorithm. To do so, we focus on the example of the car with moving obstacles in figure \ref{fig:carM}, and vary the time step $\delta$, the time horizon $T$, and the PDHG parameters $\sigma, \tau, \kappa.$ Once again, the iteration counts are averaged over 50 trials, and the average is rounded to the nearest integer. The results are in figure \ref{fig:tables}. In each table in the figure, aside from the parameter being highlighted, all other parameters are held at their baseline values of $\delta = 0.1, T = 6, \sigma = \tau = 0.5, \kappa = 1$. 

\begin{figure}[h!]
\centering \RV{
\begin{tabular}{| l l |}
\hline
$\delta$ & Iters. \\
\hline
0.1 & 1748\\
0.05 & 1777\\
0.025& 2578\\
0.0125 & 4317\\
\hline
\end{tabular} \hspace{0.5in} \begin{tabular}{| l l |}
\hline
$T$ & Iters. \\
\hline
4 & $100000+$\\
 6 & 1748\\
 10 &  864\\
14& 874\\
\hline
\end{tabular} 
 
\vspace{0.5cm}

\begin{tabular}{| l l |}
\hline
$\sigma$ & Iters. \\
\hline
 0.5 & 1748\\
 0.25 &  2340\\
0.125& 3405\\
0.0625 & 4382\\
\hline
\end{tabular}  \hspace{0.5in} \begin{tabular}{| l l |}
\hline
$\tau$ & Iters. \\
\hline
 0.5 & 1748\\
 0.25 &  1660\\
0.125& 1763\\
0.0625 & 1630\\
\hline
\end{tabular}
\hspace{0.5in} \begin{tabular}{| l l |}
\hline
$\kappa$ & Iters. \\
\hline
 1 & 1748\\
 0.75 &  2173\\
 0.5& 3618\\
0.25 & $85065^*$\\
\hline
\end{tabular}}
\caption{\RV{Tables detailing how the number of iterations until convergence depends on the different parameters for the example in figure \ref{fig:carM}. Each iteration count is averaged over 50 trials, and in each table all paramaters aside from the one being highlighted are held at baseline values. For $T = 4$, no trial resulted in convergence within the maximum iteration count of 100000. The asterisk for iteration count when $\kappa = 0.25$ is included because in that case, 18 of the 50 trials failed to converge within the maximum iteration count of 100000 so the value is artificially deflated.} }
\label{fig:tables}
\end{figure}

\RV{Looking at the table, there is a reasonably clear trend that decreasing $\delta$ increases the number of iterations requires. This is likely because with smaller $\delta$ there are simply more points that need to be resolved (and more randomness in the initialization which needs to be overcome). Using this same reasoning, one may expect that increasing $T$ also increases the number of iterations required. However, this does not hold. When decreasing $T$ to $4$, there is now no path which can even come close to reaching the final point. In cases like this, we consistently found that that algorithm would fail to converge within the maximum iteration count of 100000 (and showed no apparent signs of nearing convergence). When $T$ is increased, the iteration count significantly falls. Unfortunately, choosing such $T$ is still undesirable. When $T$ is chosen too large, since any path which reaches the endpoint is an optimal path, there are infinitely many optimal paths (assuming small-time local controllability). Thus the algorithm has an easier time finding an optimal path simply because there are more optimal paths. But observing the selected path, the car will usually dawdle for a bit, wasting some time before arriving at the endpoint at the time horizon $T$. Paths like this are likely not what one would wish to resolve, though in many cases, one could easily intuit the ``correct" optimal path from one of these paths by removing segments of the path which one deems as purposeless motion. One final note regarding either decreasing $\delta$ or increasing $T$ is that either of these changes increases the cost per iteration by increasing the number of time steps. So for example, when $\delta$ is halved, there is a two-fold increase in computation because there are twice as many points, but there is also an increase in computation due to the increased iteration count. Because of this, to maximize efficiency, it is important to choose $\delta$ as large as is allowable for ones particular application.} 

\RV{The PDHG parameters are bound by $\sigma\tau \le 0.25$ and $\kappa \in [0,1]$. Thus, from the baseline values, neither of $\sigma$ or $\tau$ can be increased without decreasing the other, and $\kappa$ can only be decreased. Looking at the tables, there are clear trends which indicate that decreasing $\sigma$ or $\kappa$ increases the number of iterations required for convergence. For $\tau$, there isn't an abundantly clear trend, though it is possible that decreasing $\tau$ could provide a marginal benefit when compared with the baseline value. Another possible benefit of decreasing $\tau$ is that it would allow one to increase $\sigma$ further, which could decrease the iteration count further. Using this observation, we tested further combinations of larger $\sigma$ with smaller $\tau$ (keeping other parameters at baseline) and from our simulations, the best case scenario was $\sigma = 0.75, \tau = 0.25$, with an average iteration count of $1442$ iterations before convergence.}

\section{Conclusion \& Future Work} \label{Conclusion}

%\begin{itemize}
%    \item Real-time methods
%    \item Modeling more realistic problems
%    \item Many agents
%\end{itemize}

In this paper, we developed an algorithm for optimal path-planning for curvature constrained motion which includes kinematic models for simple cars, airplanes, and submarines. Our method relied on a level-set, PDE-based formulation of the optimal path-planning problem, wherein the value function is not the optimal travel time (as in more control theoretic models), but optimal distance to the desired endpoint. This allowed us to solve the problem using new numerical methods which resolve the solutions to Hamilton-Jacobi equations via optimization problems. We discussed the ramifications of this modeling decision, and described in detail the implementation of our method. Finally, we demonstrated our method on some synthetic examples. In these examples, we are able to resolve optimal trajectories for cars, planes, and submarines in a matter of seconds. \RV{This allows one to maintain the PDE and control based methodology (and thus maintain interpretability) without sacrificing efficiency or scalability.}

\RV{An immediate avenue of future work is to look into manners in which convergence of Algorithm 1 could be accelerated using more sophistocated optimization methods, Nesterov-type acceleration, or better approximation of the proximal operator. Here basic gradient descent was chosen for ease of implementation and exposition, but this could very likely be improved to make the algorithm even faster and more robust.} Next, while this represents a step toward real-time PDE-based path-planning algorithms, there are still difficulties to overcome. As presented, our method can compute optimal paths in a fully-known environment. In realistic scenarios, the exact configuration of obstacles is likely unknown and obstacles may move unpredictably. \RV{Accordingly, it would be desirable if the vehicles had only local information regarding terrain and obstacles which is updated along the trajectory so that they can react and recompute paths in semi-real-time. In doing so, one would need to abandon the hope of finding globally optimal trajectories, because only local information would be known. However, this also presents the difficulty of correctly choosing the new time-horizon upon adding new information in the Hamiltonian and recomputing the optimal path, so it may require a significant change in the modeling.} In the same vein, it could prove interesting to adapt and apply our method to problems with other realistic concerns such as energy-efficient path-planning \RV{or rider comfort via jerk bounds or control regularization.} Finally, we would like to adapt this method to a multi-agent control or many-player differential games scenario.

\section*{Acknowledgments}

The authors would like to think Robert Ferrando for many discussions regarding this work. The authors were supported in part by NSF DMS-1937229 through the Data Driven Discovery Research Training Group at the University of Arizona.

%% The Appendices part is started with the command \appendix;
%% appendix sections are then done as normal sections
%\appendix

%\section{Sample Appendix Section}
%\label{sec:sample:appendix}

%% If you have bibdatabase file and want bibtex to generate the
%% bibitems, please use
%%
 \bibliographystyle{elsarticle-num} 
 \bibliography{references}

%% else use the following coding to input the bibitems directly in the
%% TeX file.

% \begin{thebibliography}{00}

% %% \bibitem{label}
% %% Text of bibliographic item

% \bibitem{}

% \end{thebibliography}
\end{document}